\documentclass[12 pt]{amsart}

\usepackage{latexsym}
\usepackage{amssymb, amsmath}
\usepackage{amsthm}
\usepackage{geometry}
\usepackage{hyperref}
\usepackage{changebar}

\numberwithin{equation}{section}

\newtheorem{theorem}{Theorem}[section]
\newtheorem{lemma}[theorem]{Lemma}
\newtheorem{cor}[theorem]{Corollary}
\newtheorem{sublem}[theorem]{Sublemma}
\newtheorem{proposition}[theorem]{Proposition}

\newtheorem{defin}{Definition}

\newtheorem{rem}[theorem]{Remark}
\newtheorem{convention}[theorem]{Convention}

\newcommand{\W}{\mathcal{W}}
\newcommand{\K}{\mathcal{K}}
\newcommand{\Lp}{\mathcal{L}}
\newcommand{\C}{\mathcal{C}}
\newcommand{\F}{\mathcal{F}}
\newcommand{\M}{\mathcal{M}}
\newcommand{\B}{\mathcal{B}}
\newcommand{\Or}{\mathcal{O}}
\newcommand{\Si}{\mathcal{S}}

\newcommand{\N}{\mathbb{N}}
\newcommand{\V}{\mathbb{V}}

\newcommand\Id{\text{\bf Id}}

\newcommand\nn{\nonumber}

\newcommand{\pd}{\partial}
\newcommand{\vf}{\varphi}
\newcommand{\po}{\varphi_1}
\newcommand{\pt}{\varphi_2}

\newcommand{\R}{\mathbb{R}}
\newcommand{\ve}{\varepsilon}
\newcommand{\tT}{\tilde{T}}
\newcommand{\ob}{\overline{\omega}}

\begin{document}

\title[Stability in two-dimensional piecewise hyperbolic maps]{Stability of statistical properties in
two-dimensional piecewise hyperbolic maps}
\author{Mark F. Demers \and Carlangelo Liverani}
\address{Mark F. Demers, Department of Mathematics, Fairfield University, Fairfield CT 06824, USA}
\email{mdemers@mail.fairfield.edu}
\address{Carlangelo Liverani, Dipartimento di Matematica, Universit\`a
di Roma {\sl Tor Vergata}, Via della Ricerca Scientifica, 00133 Roma,
Italy} 
\email{liverani@mat.uniroma2.it}
\date{\today}
\begin{abstract} 
We investigate the statistical properties of a piecewise
smooth dynamical system by studying directly the action of the transfer
operator on appropriate spaces of distributions. We accomplish such a
program in the case of two-dimensional maps with uniformly bounded
second derivative.  For the class of
systems at hand, we obtain a complete description of the SRB
measures, their statistical properties and their stability with respect
to many types of perturbations, including deterministic and random
perturbations and holes.
\end{abstract}

\thanks{We would like to thank the Institut Henri Poincar\'{e} where 
part of this work was done (during the trimester {\em Time at Work}). 
Also we enjoyed
partial support from M.I.U.R. (Cofin 05-06 PRIN 2004028108).
M.D. was partially supported by 
NSF VIGRE Grant DMS-0135290 and by the School of Mathematics of the Georgia Institute
of Technology.  
Finally, C.L. would like to
warmly thank G. Keller with whom, several years ago, he had uncountably many
discussions on these types of problems. Although we were not able to solve
the problem at the time, as the technology was not ripe yet,
the ground work we did has been precious for the present work.}


\maketitle

\section{Introduction}
\label{intro}

In recent years, many works have sought to establish in the hyperbolic
setting the functional analytic approach developed for one-dimensional piecewise expanding maps.

This strategy avoids completely any attempt to code the system and
studies directly the transfer operator on an appropriate Banach space
(in the expanding case, the functions of bounded variation). Roughly speaking, the approach is
to first obtain {\em a priori} control on the smoothing properties of the
transfer operator \cite{LY}, then infer from those that the transfer
operator is quasi-compact and that its peripheral spectrum
provides abundant information about the statistical properties of
the system \cite{Kel1}, and finally show that such a picture is stable
for a large class of perturbations \cite{BalYo, KL}. See \cite{Bal} for
a detailed explanation of the above ideas and complete references and
\cite{Li2} for an apology. 

Such a point of view was successfully extended to multidimensional
expanding maps \cite{saussol, buzzi, tsujii 1, tsujii 2, buzzi keller}, but 
its application to the
hyperbolic setting has been lacking until recently.  Notwithstanding
some partial successes \cite{Bakh, Li1,Rugh1,Rugh2,Rugh3}, the first
paper in which the above approach was systematically implemented in
all its aspects was \cite{BKL}, in which the authors studied Anosov diffeomorphisms. 
Such results have subsequently been dramatically improved in a series of
papers \cite{liverani gouezel, Ba1, BaT, Li3} of which certainly we have not seen
the end.

In spite of the fact that in one dimension the
approach was developed to overcome the problem of discontinuities, the
case of piecewise hyperbolic systems has eluded attempts to treat it
along such lines (with the partial exception of \cite{Li1}). 
Consequently, as far as hyperbolic systems with discontinuities are
concerned, the only available approaches are \cite{pesin} and
\cite{young} (and the generalizations by Chernov \cite{chernov 0} and
Chernov, Dolgopyat \cite{chernov dolgopyat} of Young's approach
in the case of billiards; see \cite{chernov young} for a review). 
Such approaches require a very deep
preliminary understanding of the regularity properties of the
invariant foliations and are not well-suited
to the study of perturbations of the systems under consideration. 

The present paper makes a first step in overcoming the difficulty of discontinuities 
by showing that
in the two-dimensional case the functional analytic approach
can be carried out successfully.  We are confident that this approach can
be extended to billiards and higher dimensional systems.

The outline of the paper is as follows. 
In Section~\ref{compactness}, we describe functional spaces on which we
establish the quasi-compactness of the transfer operator in Section~\ref{inequalities}. 
This is the
key result of the paper from which all the rest follow.

In Section~\ref{spectrum}, we show that there exists a precise relation
between the spectral picture of the transfer operator 
and the statistical properties of the
system. More precisely, the peripheral spectrum corresponds to the
ergodic decomposition with respect to the physical measures, so a complete
description of the SRB measures for the system is obtained.  

The rest of the spectrum is connected to the finer statistical properties such as
the decay of correlations, which is proven to be exponential for mixing
systems, the Central Limit Theorem, the power spectrum and the Ruelle resonances.
Although the decay of correlations and CLT are already known for systems
with a slightly more restrictive class of singularities
(see \cite{young}), the current approach presents a unified framework for these results
and adds to them a detailed understanding of the power spectrum and Ruelle resonances
not previously available.
(See \cite{ruelle 1, ruelle 2, parry poll 1, parry poll 2} 
for a discussion of Ruelle resonances in Axiom A systems.)

In addition, we answer questions concerning the stability of these statistical properties 
with respect to both deterministic and
random perturbations, as well as those obtained by introducing small holes into the system.
We prove that the
stability is of a very strong nature: all the statistical properties,
from the invariant measures to the rate of decay of correlations to the Ruelle resonances, vary
continuously with the perturbation.  The proofs of these results 
are contained in Section~\ref{perturbations}.

Contrary to \cite{liverani gouezel}, the spaces introduced here do not allow an
extensive study of the influence of the smoothness of the system on
its statistical properties. 
This may depend on the class of systems under investigation: contrary to the smooth case in which the 
degree of smoothness determines the size of the essential spectrum, it is conceivable that there is no 
difference between piecewise $\C^2$ and piecewise $\C^r$ systems. On the other hand there may be a 
difference that is not captured by our spaces. 

Finally, note that the paper tackles the problem, left open in \cite{liverani gouezel}, of how
to define spaces with {\em H\"older regularity} in the unstable direction.

\begin{rem}
A remarkable aspect of the present approach is that it bypasses
completely the detailed, and extremely laborious, study of the
smoothness properties of the invariant foliations, their
holonomies and the local ergodicity theorems
(albeit restricted to the uniformly hyperbolic case).
Accordingly, it provides an extremely direct way to obtain very strong results, as testified by the
length of the present, essentially self-contained, paper.
\end{rem}

\begin{convention}
\label{convention: C}
In this paper we will use $C$ to denote a generic constant depending
only on the dynamical systems $(\M,T)$, while $C_{a,b,c\dots}$ will
depend only on $(\M,T)$ and the parameters $a,b,c,\dots$. Accordingly,
the actual value of such constants may vary from one occurrence to the
next.
\end{convention}


\section{Setting, Definitions and Results}
\label{setting}
  
Let $\M$ be a compact two-dimensional Riemannian manifold, possibly with boundary 
and not necessarily connected, and let 
$T:\M \circlearrowleft$ be a piecewise uniformly hyperbolic map in the
following sense.
We assume that there exist a finite number of pairwise disjoint open regions $\{ \M^+_i \}$ such that
$\cup_i \overline{\M^+_i} = \M$ and the boundaries of $\M^+_i$ are piecewise $\C^1$
curves of finite length.  We define $\M_i^- = T(\M^+_i)$ and require that $\cup_i \overline{\M^-_i} = \M$.
We refer to the sets $\Si^\pm := \M \backslash \cup_i \M_i^\pm$ as
the singularity sets for $T$ and $T^{-1}$ respectively.
We assume that $T\in\operatorname{Diff}^2(\M \backslash \Si^+,\M
\backslash \Si^-)$ and that on each $\M^+_i$, $T$ has a $\C^2$ extension
to $\overline{\M_i^+}$.

On each $\M_i$, $T$ is uniformly  hyperbolic: i.e., it admits two continuous $DT$-strictly-invariant 
families of cones $C^s$ and $C^u$ defined on all of  
$\M\backslash (\Si^+ \cup \partial \M)$ which satisfy,\footnote{Note that the strict 
invariance of the cone field together with the smoothness properties of the map 
implies that the stable and unstable directions are well-defined for each point
whose trajectory does not meet a singularity line.}   
\begin{equation}
\label{eq:exp def}
\begin{split}
& \lambda \; := \; \inf_{x\in \M\backslash \Si^+} \; \inf_{v \in C^u}
\frac{\| DTv \|}{\|v\|} \; > \; 1  , \\ 
&   \mu \; := \; \inf_{x\in \M\backslash \Si^+} \; \inf_{v \in C^s}
\frac{\| DTv \|}{\|v\|} \; < \; 1  , \\
&   \mu_+^{-1} \; := \; \inf_{x\in \M\backslash \Si^-} \; \inf_{v \in C^s}
\frac{\| DT^{-1}v \|}{\|v\|} \; > \; 1  .
\end{split}
\end{equation}
In Section~\ref{leaves}, we define narrow cones with the same names
and refer to them as the stable and unstable cones of $T$  
respectively.  We assume that the tangent vectors to the
singularity curves in $\Si^-$ are bounded away from $C^s$. 
Note that this class of maps is similar to that studied in
\cite{young, pesin}; see also \cite{liverani wojtkowski} for the symplectic case.

\begin{rem}
\label{rem:alt cond}
We can replace the condition that the singularity curves be transverse to $C^s$ by the
more general assumption (H1) of Section~\ref{hole results} 
(replacing $\partial H$ with $\Si^-$), 
thus allowing singularities which are in
places tangent to the stable direction.
The estimates of Section~\ref{holes} imply that Proposition~\ref{prop:lasotayorke} and 
Theorem~\ref{thm:invariant measures} hold with this weaker condition on the singularities 
of $T$ as
long as we choose $\beta \leq \alpha/2$ in the definition of the strong unstable norm 
\eqref{eq:s-unstable}.
We do not do this, however, since this restriction on $\beta$ makes less optimal our estimates
on the essential spectral radius (see Remark~\ref{rem:optimal}). 
\end{rem}

\begin{rem}
Although the class of maps $T$ which we consider does not contain billiards,
it does contain piecewise toral automorphisms 
and a broad class of piecewise hyperbolic
nonlinear maps with bounded derivative.
\end{rem}

Denote by $\Si_n^-$ the set of singularity curves for $T^{-n}$ and by $\Si_n^+$ the 
set of singularity 
curves for $T^n$.
Let $M(n)$ denote the maximum number of singularity curves
in $\Si_n^-$ which intersect at a single point.  We make the following
assumption regarding the  singularities of $T$.

\bigskip
\noindent
\parbox{.1 \textwidth}{(P1)}
\parbox[t]{.8 \textwidth}{There exist $\alpha_0 > 0$ and an integer
$n_0>0$, such that $\lambda \mu^{\alpha_0} > 1$  
and $(\lambda \mu^{\alpha_0})^{n_0} > M(n_0)$.    } 

\bigskip
\noindent
Condition (P1) can always be satisfied if $M(n)$ has polynomial growth (as is the case
with billiards); however, since (P1) is required only for some fixed $n_0$, it is
not necessary to control $M(n)$ for all $n$ in order to verify the condition.

\begin{rem}
\label{rem:iterate p1}
If property (P1) holds for $\alpha_0$, then it holds for all
$0 < \alpha < \alpha_0$ with the same $n_0$. 
Notice also that $M(kn_0) \leq M(n_0)^k$ which implies that the
inequality in (P1) can be iterated to make  
$(\lambda \mu^{\alpha_0})^{-kn_0} M(kn_0)$
arbitrarily small once (P1) is satisfied for some $n_0$.
\end{rem}

In Section~\ref{leaves} we will define a set of admissible leaves
$\Sigma$, close to the stable direction, on which  we will define our
norms.  For a leaf $W \in \Sigma$, let $L_n$ denote the number of
smooth connected components of $T^{-n}W$.  For a fixed $N$, by
shrinking the maximum length $2\delta$ of leaves in $\Sigma$, we can
require that $L_N \leq M(N)+1$. This implies that choosing $N=kn_0$,
we can make $(\lambda \mu^{\alpha_0})^{-N}L_N$ arbitrarily small. 

\begin{convention}
\label{convention: n_0=1}
In what follows, we will assume that $n_0 = 1$.  If this is not the
case, we may always consider a higher  iterate of $T$ for which this
is so by assumption (P1).  We refer to $L_1$ as $L$ and choose $\delta$ small 
enough that $L \lambda^{-1} \mu^{-\alpha_0}=:\rho<1$.
\end{convention}

We write $D^s$ to denote differentiation in the stable
direction  and note that
this direction is well-defined outside the set $\cup_{n \geq 0} \Si_n^+$
due to the uniform hyperbolicity of $T$.

For an admissible leaf $W \in \Sigma$, we will
denote by $m$ the (unnormalized) Riemannian volume on $W$ and by
$d(\cdot, \cdot)$ the distance along the leaf. 
We will often abbreviate $m(W)$ by $|W|$.


\subsection{Transfer Operator}
\label{transfer}
The basic object of study in the present paper is the so-called {\em transfer
operator} $\Lp$. Clearly, to make sense of an operator it is necessary to
specify on which space it acts. In fact, the search for a good space is
the main point of the present paper. 

In the smooth case \cite{liverani gouezel}, it is convenient to define
the transfer operator acting on the space of distributions which
turns out to contain all the 
relevant spaces. In this manner one can obtain all the relevant
operators as restrictions of the original one.

In the present case it is not clear if there exists an appropriate
ambient space.\footnote{Clearly the space of distributions will not do
since if $\vf$ is smooth, $\vf\circ T$ may not be.} We bypass this
problem by defining the operator as acting between two scales of spaces.

For each $n\in\N$, let $\mathcal K_n$ be the set of connected components of
$\M\setminus \Si^+_n$.
Recall that $\C^1(\overline{K}, \R)$ is the set of functions $\varphi \in \C^1(\overset{\circ}{K}, \R)$
which have a $\C^1$ extension in a neighborhood of $\overline{K}$.
Let $\C^1_{\Si_n^+}:=\{\vf\in L^\infty(\M)\;:\;
\vf\in\C^1(\overline{K},\R)\;\forall K\in\mathcal K_n\}$.\footnote{The space
$\C^1_{\Si_n^+}$  is a Banach space when equipped with the norm
$\sup\limits_{K\in\mathcal K_n}|\vf|_{\C^1(\overset{\circ}{K})}$.}
If $h\in(\C^1_{\Si_n^+})'$, is an element of the dual of $\C^1_{\Si_n^+}$,
then
$\Lp:(\C^1_{\Si_{n}^+})'\to (\C^1_{\Si_{n-1}^+})'$ acts on $h$ by
\[
\Lp h(\vf) = h(\vf \circ T) \quad \forall \vf\in \C^1_{\Si_{n-1}^+}.
\]

The above definition shows how the transfer operator acts on an
abstract space of distributions, but often we will be concerned with its
action on more concrete objects. Notice that since the sets
$\Si^+_n$ are all of zero Lebesgue (Riemannian) measure, each signed measure
absolutely continuous with respect to Lebesgue yields an element of
$(\C^1_{\Si_n^+})'$. 

\begin{rem}
In what follows, we will identify a measure $h$ that is
absolutely continuous with respect to Lebesgue with its density,
which we will insist on calling $h$.  Accordingly,
\[
h(\vf)=\int_\M h\vf\, dm
\]
where $m$ denotes Lebesgue measure on $\M$. Hence the space of  measures
absolutely continuous with respect to Lebesgue is canonically identified
with $L^1(\M,\R,m)$.
\end{rem}

With the above convention, $L^1(\M)\subset (\C^1_{\Si_n^+})'$ for each
$n\in\N$. 
One can then restrict $\Lp$ to $L^1$ and a simple
computation shows that\footnote{Given a square matrix $A$, by $|A|$ we
mean $|\det(A)|$.}
\[
\Lp^n h = h \circ T^{-n} \; |DT^n(T^{-n})|^{-1}
\]
for any $n \geq 0$ and any $h \in L^1(\M)$.\footnote{Often the above
is taken as the definition of the transfer operator, yet as will
become clear in the following, $L^1$ is both too small and too large a space to be useful.}


\subsection{Definition of the Norms}
\label{norms}

We will define the required Banach spaces by closing $\C^1$ with respect to
suitable norms.

The norms are defined via a set of {\em admissible leaves}
$\Sigma$. Such leaves are essentially smooth curves roughly in
the stable direction, their length is smaller than some $\delta$ and
among them is defined a notion of distance $d_\Sigma$.  Also, a notion of
distance $d_q$ is defined among functions supported on such leaves.
They are defined precisely in Section~\ref{leaves}.

For $W \in \Sigma$ and $0 \leq \alpha,q \leq 1$, denote by
$\C^\alpha(W, \mathbb{C})$ the set of continuous complex-valued functions on 
$W$ with H\"{o}lder exponent $\alpha$.
Define the following norms
\[
|\vf|_{W,\alpha,q}:=|W|^\alpha\cdot|\vf|_{\C^q(W,\mathbb{C})} .
\]

Given a function $h \in \C^1(\M, \mathbb{C})$, define the \emph{weak norm} 
of $h$ by
\begin{equation}
\label{eq:weak}
|h|_w:=\sup_{W\in\Sigma}\sup_{\substack{\vf\in\C^1(W,\mathbb{C})\\
|\vf|_{\C^1(W)}\leq 1}}\int_W h \vf \; dm .
\end{equation}
Choose $\alpha$, $\beta$, $q < 1$ such that 
$0 < \beta \leq \alpha \leq 1-q \leq
\alpha_0$.\footnote{Such inequalities are irrelevant
for the definition of the spaces, but we introduce them here because they will
be needed for the various Lasota-Yorke estimates of Section~\ref{inequalities}.}
We define the \emph{strong stable norm} as
\begin{equation}
\label{eq:s-stable}
\|h\|_s:=\sup_{W\in\Sigma}\sup_{\substack{\vf\in\C^1(W,\mathbb{C})\\
|\vf|_{W,\alpha,q}\leq 1}}\int_W h \vf \; dm
\end{equation}
and the \emph{strong unstable norm} as
\begin{equation}
\label{eq:s-unstable}
\|h\|_u:=\sup_{\varepsilon \leq \varepsilon_0} \; \sup_{\substack{W_1,
W_2 \in \Sigma \\
d_\Sigma (W_1,W_2)\leq \varepsilon}}\;
\sup_{\substack{|\vf_i|_{\C^1(W_i,\mathbb{C})}\leq 1\\ d_q(\vf_1,\vf_2)\leq \ve}} \;
\frac{1}{\varepsilon^\beta} \left| \int_{W_1} h
\vf_1 \; dm - \int_{W_2} h \vf_2 \; dm \right| 
\end{equation}
where $\ve_0$ will be chosen later.
We then define the \emph{strong norm} of $h$ by
\begin{equation}
\label{eq:strong}
\|h\| = \|h\|_s + b \|h\|_u
\end{equation}
where $b$ is a small constant chosen in Section~\ref{inequalities}.

We define $\B$ to be the completion of $\C^1(\M)$ in the strong norm and
$\B_w$ to be the completion of $\C^1(\M)$ in the weak norm.

Finally, let
\begin{equation}
\label{eq:def-Dn}
D_n:=\delta^{\alpha - 1}\sup_{0\leq k\leq n}\;\sup_{W\in\Sigma}|W|^{-\alpha}\int_W|DT^{-k}|dm
\end{equation}
and set $D_* = \limsup_{n\to \infty} D_n^{1/n}$.


\subsection{Statement of Results}
\label{results}

The first result gives a more concrete description of the above
abstract spaces.

\begin{lemma}
For each $n \geq 0$, $\B\subset \B_w\subset (\C^1_{\Si^+_n})'$.
\end{lemma}
\begin{proof}
This is an immediate consequence of Lemma~\ref{lem:distribution}
and the fact that $|\cdot|_w\leq \|\cdot\|$.
\end{proof}

In addition, the transfer operator is well-defined on the spaces
$\B,\;\B_w$. In fact, the following more precise result is proven in
Section \ref{inequalities}. 

\begin{proposition}
\label{prop:lasotayorke}
There exists $\delta_0>0$ 
such that for all $h \in \B$, $\delta\leq\delta_0$ and $n\geq 0$,
\begin{eqnarray}
\label{eq:weak norm}  
| \Lp^n h |_w   & \leq & C D_n |h|_w \; , \\     
\label{eq:stable norm} 
\| \Lp^nh\|_s & \leq & C\max\{\rho,\mu_+^q\}^nD_n\|h\|_s + C_\delta
D_n|h|_w\; ,  \\      
\label{eq:unstable norm}    
\| \Lp^nh\|_u & \leq & C \lambda^{-\beta n} D_n \|h\|_u + C(D_n+L_n \lambda^{-n} \mu^{-\alpha n}) \|h\|_s .
\end{eqnarray}
\end{proposition}

If we choose $1>\tau>\max\{\lambda^{-\beta}, \rho,\mu_+^q\}$, then there exists
$N \geq 0$ such that
\begin{equation}
\label{eq:almost-ly}
\begin{split}
\| \Lp^N h\| & =  \| \Lp^N h\|_s + b \| \Lp^N h\|_u \\   
          &  \leq  \frac{\tau^N D_N}2 \|h\|_s + C_\delta D_N |h|_w + b \tau^N D_N \|h\|_u 
          + b C(D_N+L_N \lambda^{-N} \mu^{-\alpha N}) \|h\|_s \\
          &  \leq  \tau^N D_N \|h\| + C_\delta D_N |h|_w
\end{split}
\end{equation}
provided $b$ is chosen small enough with respect to $N$.
The above represents the
traditional Lasota-Yorke inequality once we show the $D_n$ are
bounded. Probably a direct argument could prove this fact, yet we
find it easier to prove using a functional analytic argument. 

The final ingredient in the strategy to prove the
quasi-compactness of the operator $\Lp$ is the relative compactness
of the unit ball of $\B$ in $\B_w$.  This is proven in Lemma \ref{lem:compact}. It thus follows
by standard arguments (\cite{Bal}) that the essential spectral radius
of $\Lp$ on $\B$ is bounded by $\tau D_*$, while the estimate
for the spectral radius, contrary to the usual situation, is $D_*$
which, in general, could be
larger than one. Nevertheless, a functional analytic argument (Lemma
\ref{lem:spectra}) shows that the spectral radius is one. As a
consequence we know, {\it a posteriori}, that $D_*=1$ 
and this together with Lemma~\ref{lem:D*} implies that the $D_n$
are bounded (see Remark~\ref{rem:D*}). 

Our next results characterize the set of invariant measures in $\B$ and some of the 
statistical properties of $T$.  Recall
that an invariant  probability measure $\mu$ is called a {\em physical
measure} if there exists a 
positive Lebesgue measure invariant set $B_\mu$, with $\mu(B_\mu)=1$, 
such that, for each continuous
function $f$,
\[
\lim_{n\to\infty}\frac 1n\sum_{i=0}^{n-1}f(T^ix)=\mu(f)\quad \forall x
\in B_\mu.
\]
Let $\Pi_\theta$ be the eigenprojector on $\V_\theta$, the eigenspace of $\Lp$ corresponding to eigenvalue 
$e^{2\pi i\theta}$, and set $\V:=\oplus_{\theta}\V_\theta$.
The following theorem is proved by the lemmas of Section~\ref{spectrum}.

\begin{theorem}
\label{thm:invariant measures}
The peripheral spectrum of $\Lp$ on $\B$ consists of finitely many
cyclic groups. 
The maps $\{T^n\}_{n\in\N}$ admit only finitely many physical 
probability measures, they 
form a basis for $\V$ and the cycles correspond to the cyclic groups. 
In addition,
\begin{enumerate}
  \item If $\mu \in \V_0$ and $\Si^\pm_{n,\epsilon}$ is an $\epsilon$-neighborhood of $\Si^\pm_n$, 
then $\mu(\Si^\pm_{n,\epsilon}) \leq C_n \epsilon^\alpha$ for all $n \in \mathbb{N}$.  In particular,
$\mu(\Si^\pm_n)=0$.
  \item Each element in $\V$ is a signed measure absolutely continuous with respect to the probability measure  
$\bar{\mu}  := \lim_{n\to\infty} \frac{1}{n} \sum_{i=0}^{n-1} \Lp^i 1$. In particular, all the physical measures are   
absolutely continuous with respect to $\bar{\mu}$.
  \item The supports of the physical measures correspond to the ergodic decomposition
with respect to Lebesgue.
  \item For all $f \in \C^0(\M, \mathbb{R})$, the limit 
    $f^+(x) := \lim_{n\to \infty} \frac{1}{n} \sum_{i=0}^{n-1}
f\circ T^i(x)$ exists for $m$-almost-every 
    $x$ and takes on only finitely many different values.  If $\bar{\mu}$ is ergodic, then
    $f^+(x) = \int f d\bar{\mu}$ for $m$-almost-every $x$.
  \item If $(T,\bar\mu)$ is ergodic, then 1 is a simple eigenvalue. If
$(T^n, \bar{\mu})$ is ergodic for all $n\in\N$, then one is the only
eigenvalue of modulus one, $(T,\bar{\mu})$ is mixing and exhibits
exponential decay of correlations for H\"older
observables, and the Central Limit Theorem holds.
\item More generally, the Fourier transform of the correlation
function (sometimes called the {\em power spectrum}) admits a meromorphic 
extension in the annulus $\{z\in{\mathbb C}\;;\; \tau<|z|<\tau^{-1}\}$ and the 
poles (Ruelle resonances)
correspond exactly to the eigenvalues of $\Lp$. 
 \end{enumerate} 
\end{theorem}

Items (1-4) and part of (5) are proved in Section~\ref{peripheral}.  
The rest is proved in 
Section~\ref{statistics}.

\begin{rem}
Although $\bar{\mu}$ is a natural measure in the sense that it is obtained by pushing forward
and averaging Lebesgue measure, it is generally not absolutely continuous with respect to
Lebesgue in the hyperbolic setting.  Typically, one expects $\bar{\mu}$ to be
singular along stable manifolds and absolutely continuous along unstable manifolds.
\end{rem}

\begin{rem}
A natural question is if all the positive elements of $\V_0$ are SRB
measures; however, the characterization of SRB measures as measures
that are absolutely continuous along unstable manifolds is a bit at odds
with our philosophy since it would require us to prove the existence
and properties
of such manifolds in the first place.  An alternative approach
is to note that the integral along a manifold lying in the unstable
cone yields an element of $\B$ (see \cite[Proposition 4.4.]{liverani gouezel}
for a similar result in that context) and therefore iterating it
(one standard manner to construct SRB measures) one converges to
the elements of $\V_0$.  With this approach one can show that
$\V_0$ corresponds exactly to the decomposition into SRB measures.
\end{rem}

\begin{rem}
Several of the above results are similar to those obtained in \cite{pesin, young} for
piecewise hyperbolic maps.  
In \cite{young},
an SRB measure $\nu$ was constructed and under the assumption that $(T^n, \nu)$ is
ergodic for all $n$, it was proven that $(T, \nu)$ satisfies the CLT
and exponential decay of correlations for H\"{o}lder observables.  In \cite{pesin},
the existence of SRB measures and the ergodic decomposition was proven. 
\end{rem}

In Section \ref{perturbations}, we prove various
perturbation results, using the framework provided by \cite{KL}.  This requires first obtaining
uniform Lasota-Yorke estimates for the perturbed operators $\Lp_\ve$.  Then, regarding these operators
as acting from $\B$ to $\B_w$, we define the norm
\[
||| \Lp ||| = \sup_{ \{ h \in \B: \|h\| \leq 1 \} } |\Lp h|_w
\]
and show that $\Lp_\ve$ and $\Lp$ are close in this norm.  The results of \cite{KL} then imply
that the spectral picture (hence the SRB measures, the rate of correlation decay, etc.) persists
and is stable as long as a spectral gap is maintained.  These results, to our knowledge, are new and are a
simple byproduct of the present approach.


\subsection{Deterministic and Smooth Random Perturbations}
\label{pert results}

We define the class of perturbations for which our results hold.  This class is analogous
to that studied in \cite{liverani gouezel}.

Fix $B_* < \infty$ such that $|D^2T| < B_*$ and
let $\Gamma_{B_*}$ be the set of maps $\tT$ that satisfy the assumptions of Section~\ref{setting}
with $|D^2\tT| \leq B_*$.   

\begin{defin}
Given two maps $T_1,\,T_2 \in \Gamma_{B_*}$ we say that they have distance $\ve$ if
their singularity curves are at distance $\ve$ and if outside an $\ve$
neighborhood of the union of their singularity curves they are $\ve$-close 
in the $\C^2$ norm.  We call this distance $\gamma(T_1, T_2)$.
\end{defin}

Choose $\ve \leq \ve_0$ and let $X_\ve$ be an $\ve$-neighborhood of $T$ in $\Gamma_{B_*}$,
\[
X_\ve = \{ \tT \in \Gamma_{B_*}: \gamma(T, \tT) < \ve \}.
\]
In general, the constants $\lambda(\tT)$, $\mu(\tT)$, $\mu_+(\tT)$ and $D_n(\tT)$ defined by 
\eqref{eq:exp def} and \eqref{eq:def-Dn} depend on the map $\tT$.  However, for
$\ve \leq \ve_0$, we may choose constants $\lambda$, $\mu$, $\mu_+$ and $D_n$
such that $1 < \lambda \leq \lambda(\tT)$, $1 > \mu \geq \mu(\tT)$,
$1 > \mu_+ \geq \mu(\tT)$ and $D_n \geq D_n(\tT)$ for all $\tT \in X_\ve$.
These are the constants we shall use in the estimates of Section~\ref{perturbations}
which enable us to obtain uniform Lasota-Yorke type inequalities for the maps in $X_\ve$.

Let $\nu$ be a probability measure on a probability space $\Omega$ and let
$g:\Omega \times \M \to \mathbb{R}^+$ be a measurable function satisfying:
\begin{itemize}
  \item[(i)] $g(\omega, \cdot) \in \C^1(\M, \mathbb{R}^+)$ for each $\omega \in \Omega$;
  \item[(ii)] $\int_\Omega g(\omega, x) d\nu(\omega) = 1$ for each $x \in \M$;
  \item[(iii)] $g(\omega,x) \geq a > 0$ and $|g(\omega, \cdot)|_{\C^1(\M)} \leq A < \infty$.
\end{itemize}
If we associate to each $\omega \in \Omega$ a map $T_\omega \in X_\ve$, this defines a random walk on 
$\M$ in a natural way.  Starting at $x$, we choose $T_\omega$ according to the distribution
$g(\omega,x)d\nu(\omega)$.  We apply $T_\omega$ to $x$ and repeat this process starting at $T_\omega x$.
We say the process has size $\Delta(\nu,g) \leq \ve$.

\begin{rem}
If $\nu$ is a Dirac measure centered at $\omega_0$, this process corresponds to the deterministic
perturbation $T_{\omega_0}$ of $T$.  Thus this setting encompasses a large class of random and
deterministic perturbations of $T$.
\end{rem}

The transfer operator $\Lp_{\nu,g}$ associated with the random process governs the evolution
of densities by
\[
\Lp_{\nu,g} h(x) = \int_\Omega \Lp_{T_\omega}h(x) \, g(\omega, T_\omega^{-1}x ) \; d\nu(\omega)
\]
where $\Lp_{T_\omega}$ is the transfer operator associated with $T_\omega$.

Lemmas \ref{lem:perturbation}, \ref{lem:random perturbation} and \ref{lem:uniform estimates}
prove the two steps required in order to apply \cite{KL} to the above
class of perturbations.  We need some more notation before stating the theorem fully.

Choose $\sigma \in (\max\{\lambda^{-\beta}, \rho, \mu_+^q\}, 1)$ and denote by sp$(\Lp)$ the spectrum
of $\Lp$ on $\B$.  Since sp$(\Lp) \cap \{ z \in \mathbb{C}: |z| \geq \sigma \}$ consists
of a finite number of eigenvalues $\varrho_1, \ldots, \varrho_k$ of finite multiplicity, we may assume
that sp$(\Lp) \cap \{ z \in \mathbb{C}: |z|=\sigma \} = \emptyset$.  Hence there exists $t_*>0$
such that $|\varrho_i - \varrho_j| > t_*$ for $i \neq j$ and dist(sp$(\Lp),\{|z|=\sigma\}) > t_*$. 

Finally, for $t \leq t_*$, define the spectral projections
\[ \begin{split} & \Pi_{\nu,g}^{(j)} := 
      \frac{1}{2\pi i} \int_{|z-\varrho_j| = t} (z - \Lp_{\nu,g})^{-1} dz \qquad \mbox{and} \\
      & \Pi_{\nu,g}^{(\sigma)} := 
      \frac{1}{2\pi i} \int_{|z|=\sigma} (z - \Lp_{\nu,g})^{-1} dz.
    \end{split} \]
We denote by $\Pi_0^{(j)}$ and $\Pi_0^{(\sigma)}$ the corresponding spectral projections for
the unperturbed operator $\Lp$.

\begin{theorem}
\label{thm:pert}
For each $t \leq t_*$ and $\eta < 1- \frac{\log \sigma}{\log \max\{\lambda^{-\beta}, \rho, \mu_+^q\} }$,
there exists $\ve_1 > 0$ such that for any perturbation $(\nu,g)$ of $T$ satisfying
$\Delta(\nu,g) < \ve_1$, the spectral projections $\Pi_0^{(j)}$, $\Pi_0^{(\sigma)}$,  
$\Pi_{\nu,g}^{(j)}$ and $\Pi_{\nu,g}^{(\sigma)}$ are well-defined and satisfy
\begin{enumerate}
  \item $||| \Pi_{\nu,g}^{(j)} - \Pi_0^{(j)}||| \leq C \Delta(\nu,g)^\eta$ and
    $||| \Pi_{\nu,g}^{(\sigma)} - \Pi_0^{(\sigma)}||| \leq C \Delta(\nu,g)^\eta$;
  \item $\mbox{rank}(\Pi_{\nu,g}^{(j)}) = \mbox{rank}(\Pi_0^{(j)})$ for each $j$;
  \item $\|\Lp^n_{\nu,g}\Pi_{\nu,g}^{(\sigma)} \| \leq C \sigma^n$ for all $n \geq 0$.
\end{enumerate}
\end{theorem}

In view of the previous discussion on the meaning of the spectral
data, Theorem \ref{thm:pert} implies that the statistical properties
(invariant measures, rates of decay of correlations, variance of the
CLT, etc.) are stable under the above class of perturbations.

\begin{rem}
It is possible to obtain a constructive bound on $\ve_1$ by estimating $\tau$ 
and using the bounds provided by \cite{KL}.
\end{rem}


\subsection{Hyperbolic Systems with Holes}
\label{hole results}

Another interesting class of perturbations is the one obtained by
opening small holes in the system, thus making it an open system from which particles or
mass can escape.  In such systems, we keep track of the iterates
of points as long as they do not enter the holes.
 
Let $H \subset \M$ be an open set which we call the hole and define $\M^0 = \M \backslash H$.
Let $\M^n = \cap_{i=0}^n T^i\M^0$ be the 
set of points that has not escaped by time $n$.  The map $\tT^n := T^n|\M^n$ describes the dynamics
in the presence of the hole and
the evolution of measures is described by the transfer operator
\[
\Lp^n_Hh=\Lp^n(1_{\M^n} h).
\]
Since $\tT$ is simply a restriction of $T$, the family of admissible leaves $\Sigma$ does not change. 
Let $r = \sup \{ |W| : W \subset H, W \in \Sigma \}$, i.e.\ $r$ is the largest ``diameter'' of $H$ where
length is measured along admissible leaves. 

We make the following two assumptions on the hole.

\bigskip
\noindent
\parbox{.1 \textwidth}{(H1)}
\parbox[t]{.8 \textwidth}{$H$ is comprised of a finite number of open,
connected components whose boundaries consist 
of finitely many piecewise smooth curves.  Moreover, for each smooth
component $\omega$ of $\pd H$ and any point $x \in \omega$, either
\begin{enumerate}
  \item the tangent to $\omega$ at $x$ is bounded away from $C^s(x)$, or 
  \item the curvature of $\omega$ at $x$ is greater than $B$ (in the definition of $\Xi$ from 
    Section~\ref{leaves}).
\end{enumerate}  } 

\bigskip
For any $W \in \Sigma$,
let $P_n$ be the maximum number of connected components of $T^{-n}W \cap \cup_{i=0}^n T^{-i}H$.

\bigskip
\noindent
\parbox{.1 \textwidth}{(H2)}
\parbox[t]{.8 \textwidth}{There exists an integer $n_1>0$, such that
$(\lambda \mu^{\alpha_0})^{n_1} > P_{n_1}$.    }  

\bigskip
Notice that we can iterate the inequality in (H2) by controlling $\delta$.  For a fixed $N = kn_1$, we can
choose $\delta$ so that $P_{kn_1} \leq P_{n_1}^k$.  Thus we can make $P_N (\lambda \mu^{\alpha_0})^{-N}$
as small as we like.

\begin{convention}
\label{convention: n_1=1}
We will assume that $n_1 = 1$.  If this is not the case, we can always consider
a higher iterate of $T$ for which this is true once (H2) is satisfied.
We refer to $P_1$ as simply $P$ and assume that 
$\lambda^{-1} \mu^{-\alpha_0} (L+P) <1$.
\end{convention}

The observations following (H2) and (P1) imply that we can control
$\lambda^{-n} \mu^{-\alpha n}(L_n + P_n)$ which is precisely what we need in order to prove the
Lasota-Yorke inequalities for $\Lp_H$.  

\begin{rem}
It is fairly easy to have holes that satisfy our
assumptions: for example holes with boundaries transverse to the
stable cones, convex holes with boundaries with curvature larger than
$B$ or some appropriate mixture of the two.  In the case of
convex holes, $P =1$.
\end{rem}

\begin{rem}
We do not distinguish between pieces of $\tT^{-n}W$ created by intersections with
the hole and those created by the singularities of $T$.  This is clear
in the estimates of Sections~\ref{inequalities} 
and \ref{holes} and justifies Remark~\ref{rem:alt cond} that all the theorems
of Section~\ref{results} hold with the weaker conditions (H1) and (H2) on $\Si^-$ as
long as we choose $\beta \leq \alpha/2$.  
This restriction on $\beta$ stems from the observation that for strictly convex holes (singularities),
two curves in $\Sigma$ which are $\ve$-close to one another can differ in their
intersection with the hole (singularity) by a length of at most $C\ve^{1/2}$.  (This is used in
equation~\eqref{eq:unstable hole} of Section~\ref{holes}.)
\end{rem}

The spectral radius of $\Lp_H$ is typically $\vartheta <1$ when all the mass in the system
eventually escapes.  The analogous notion to an invariant measure in this setting is that of a
\emph{conditionally invariant measure}.  
For any Borel measure $\mu$, define $\tT_*\mu(A) = \mu(\tT^{-1}A)$ for any Borel set
$A \subset \M$. A probability measure $\mu$
is called conditionally invariant with respect to 
$\tT$ if $\tT_*\mu = \lambda \mu$ for some $\lambda \leq 1$.  It follows that $\lambda = \mu(\M^1)$ and
that $-\log \lambda$ represents the exponential rate of escape from the system with respect to $\mu$.

In principle there can be many conditionally invariant measures with different eigenvalues; however, one can
ask if there exists a natural conditionally invariant measure which is
the forward limit of a reasonable class 
of measures under the nonlinear operator $\tT_*^n\mu/|\tT_*^n\mu|$ 
(see \cite{demers young} for a discussion of the issues
involved).  Lemma~\ref{lemma:hole close} and Proposition~\ref{prop:hole compact} place us in the setting of \cite{KL}
and allow us to assert the following theorem.

\begin{theorem}
\label{theorem:holes}
Let $H$ be a hole satisfying conditions (H1) and (H2) and choose $\beta \leq \alpha/2$ in
\eqref{eq:s-unstable}.  Then for $Pr^\alpha$ sufficiently small,
\begin{enumerate}
  \item The non-essential spectra and the relative spectral projectors of $\Lp$ and $\Lp_H$ outside the 
    disk of radius $\tau$ are close in the sense of Theorem~\ref{thm:pert}.
  \item If $T$ has a unique SRB measure, then $\tT$ admits a unique
    natural conditionally invariant measure $\mu$ which is characterized
    by $\mu = \lim_{n \to \infty}
    \tT_*^nm/|\tT_*^nm|$.  
\end{enumerate}
\end{theorem}

\begin{cor}
Suppose $T$ has a unique SRB measure $\mu_0$ and
let $H_t$ be a sequence of holes with diam$(H_t) \leq t$ satisfying (H1) and (H2)
with uniform constant $n_1$.
Let $\mu_t$ be the natural conditionally invariant measures associated with $H_t$
given by Theorem~\ref{theorem:holes}(2).  Then
$|\mu_t - \mu_0|_w \rightarrow 0$ as $t\to 0$.
\end{cor}

\begin{proof}  The convergence follows directly from the closeness of the spectral projectors
guaranteed by Theorem~\ref{theorem:holes}(1).
Note that the convergence in the $|\cdot |_w$-norm is stronger than the weak-convergence results 
typically obtained for open systems.
\end{proof}

When $T$ has a unique SRB measure, one can also associate to the conditionally invariant
measure $\mu$ a unique invariant 
measure $\nu$ for $\tT$ which is supported on
$\Omega = \cap_{n=-\infty}^{\infty} T^n\M^0$, the set of points that never escape from the system.
Define $\Pi_\vartheta$ to be the projector onto the eigenspace
associated with the spectral radius $\vartheta$. 
$\Pi_\vartheta$ admits the following characterization,
\[
\Pi_\vartheta = \lim_{n\to \infty} \vartheta^{-n} \Lp_H^n.
\]
In fact, the spectral
decomposition implies that $\Lp_Hh=\vartheta \mu \ell(h)+Rh$, where the
spectral radius of $R$ is strictly smaller than $\vartheta$ and 
\[
\ell(h) = \int \Pi_\vartheta h \, dm = \lim_{n\to \infty} \vartheta^{-n} \int_{\M^n} h \, dm. 
\]
It is then easy to see that
\[
\nu(\vf) := \ell(\vf \mu) = \lim_{n\to \infty} \vartheta^{-n} \int_{\M^n} \vf \, d\mu
\]
is the required invariant measure.

\begin{rem}
Hyperbolic systems with holes have been well-studied when the systems
in question admit a finite Markov partition (see the long series of
papers  \cite{cencova, chernov 1, chernov 2, chernov 3, chernov 4, lopes
markarian}), but these are the first results
for hyperbolic systems with discontinuities and no Markov properties.
Moreover,  it should be noted that even
if $T$ is a $\C^2$ Anosov diffeomorphism, then the present approach
yields stronger results in a much more simple, 
direct and compact way than has previously been available.
In one dimension, piecewise expanding maps with non-Markov holes have
been studied via a variety of approaches, \cite{chernov bedem},
\cite{liverani maume}, \cite{demers exp}; logistic maps with non-Markov
holes were studied in \cite{demers logistic}.
\end{rem}


\section{Banach space embeddings}
\label{compactness}

We must start with the overdue exact definition of the family of admissible leaves
$\Sigma$,
which is a set of parametrized curves in the unstable direction.


\subsection{Family of Admissible Leaves}
\label{leaves}
  
Our definitions are similar to those of \cite{liverani gouezel}.

For $\kappa$ sufficiently small, we redefine the stable cone at $x \in \M$ to be
\[
C^s(x) = \{ u + v \in T_xM : u \in E^s(x), v \perp E^s(x), \| v \| \leq \kappa \| u \| \}.
\]
An analogous expression defines $C^u(x)$.
These families of cones are invariant, that is $DT^{-1}(x)(C^s(x)) \subset C^s(T^{-1}x)$ and
$DT(x)(C^u(x)) \subset C^u(Tx)$.

For each $\M_i^+$,
we choose a finite number of coordinate charts $\{ \chi_j \}_{j=1}^K$, whose domains
$R_j$ vary depending on whether they contain a preimage of part of the boundary curves of
$\M_i^+$.  For those $\chi_j$ which map only to the interior of $\M_i^+$, we take
$R_j = (-r_j,r_j)^2$.  For those $\chi_j$ which map to a part of $\pd \M_i^+$, we take
$R_j$ to be $(-r_j, r_j)^2$ restricted to one side of a piecewise $\C^1$ curve (the preimage
of part of $\pd \M_i^+$) which we position so that it passes through the origin. 
Each $R_j$ has a centroid, $y_j$, and each $\chi_j$ satisfies
\begin{enumerate}
  \item[(1)] $D\chi_j(y_j)$ is an isometry;
  \item[(2)] $D\chi_j(y_j) \cdot (\R \times {0}) = E^s(\chi_j(y_j))$;
  \item[(3)] The $\C^2$-norm of $\chi_j$ and its inverse are bounded by $1+\kappa$;
  \item[(4)] There exists $c_j \in (\kappa, 2\kappa)$ such that the cone 
    $C_j = \{ u+v \in \R^2 : u \in \R \times \{0\},
    v \in \{0\} \times \R, \|v\| \leq c_j \|u\| \}$ has the following property:  for $x \in R_j$
    such that $\chi_j(x) \notin \Si^-$,  $D\chi_j(x)C_j \supset C^s(\chi_j(x))$ and 
    $DT^{-1}(D\chi_j(x)C_j) \subset C^s(T^{-1} \circ \chi_j(x))$;
  \item[(5)] $\M_i^+$ is covered by the sets 
    $\{ \chi_j(R_j \cap (-\frac{r_j}{2},\frac{r_j}{2})^2) \}_{j=1}^K$.
\end{enumerate}

Now choose $r_0 \leq \min_j r_j/2$; later, we may shrink $r_0$ further.  
Fix $B < \infty$ and consider the set of functions
\[ 
\Xi := \{ F \in \C^2([-r,r], \mathbb{R}) : \; r \in (0, r_0], F(0) =
0, |F|_{\C^1} \leq \kappa, |F|_{\C^2} \leq B \}. 
\]

Let $I_r = (-r,r)$.  For $x \in R_j \cap (-r_j/2,r_j/2)^2$ such that $x + (t,F(t)) \in R_j$ 
for $t \in I_r$, define $G(x,r,F)$ to be a lift of the graph of 
$F$ to $\M$: $G(x,r,F)(t):= \chi_j(x + (t,F(t)))$ 
for $t \in I_r$.  For ease of notation, we will often write 
$G_F$ for $G(x,r,F)$.
We record here for future use that $|G_F|_{\C^1} \leq (1+\kappa)^2$ and $|G_F^{-1}|_{\C^1} \leq 1+\kappa$.

Our set of \emph{admissible leaves} is then defined as follows,
\[
\Sigma := \{ W= G(x,r,F)(I_r) : x \in R_j \cap (r_j/2, r_j/2)^2, r \leq r_0, F \in \Xi \}.
\]
If necessary, we shrink $r_0$ so that $\sup_{W \in \Sigma} |W| \leq
2\delta$ where $\delta$ is the length scale  
referred to in the convention following property (P1).

We define an analogous family of approximate unstable leaves $\F^u$
which lie in the unstable cone $C^u$.   

For any two leaves $W_1(\chi_{i_1}, x_1,r_1,F_1)$ and $W_2(\chi_{i_2},
x_2,r_2,F_2)$ with $r_1 \leq r_2$,
we define the distance between them to be\footnote{The reader can
check that the triangle inequality holds in $\Sigma$.}
\[
d_\Sigma (W_1,W_2) = \eta(i_1, i_2) + |x_1 - x_2| + |r_1 - r_2| +2^{-1} B^{-1}|F_1 -
F_2|_{\C^1(I_{r_1})}
\]
where $\eta(i,j) = 0$ if $i=j$ and $\eta(i,j) = \infty$ otherwise,
i.e., we can only compare leaves which are  
mapped under the same chart.

Given two functions
$\vf_i\in\C^q(W_i,\mathbb{C})$, we define the distance between
$\vf_1$, $\vf_2$ as
\[
d_q(\po,\pt) =|\po\circ G_{F_1}-\pt\circ G_{F_2}|_{\C^q(I_{r_1},\mathbb{C})}. 
\]


\subsection{Some Technical Facts}
\label{sub:technical}
To understand the structure of the spaces $\B_w$ and $\B$ it is necessary
to prove two preliminary results that will be needed in many
other arguments throughout the paper. In particular we need some
understanding of the properties of $T^{-n}W$ for $W\in \Sigma$.
We use the distortion bounds of Appendix A throughout Sections~\ref{compactness}
and \ref{inequalities}.

Let $\W_0 = \{ W \}\subset \Sigma$ and suppose we have 
defined $\W_{n-1}\subset \Sigma$.  If $W' \in \W_{n-1}$ contains
any singularity points of $T^{-1}$, then  
$T^{-1}W'$ is partitioned into at most $L$ pieces $W'_i$, so that $T$
is smooth on each $W'_i$. 
Next, if one of the components of $T^{-1}W'$ 
has length greater than $2\delta$, it is partitioned further into
pieces of length between $\delta$ and $2\delta$.
We define $\W_n$ to be the collection of all pieces $W_i \subset T^{-n}W$ obtained
in this way.  
It is a standard result of hyperbolic theory that each $W_i$ is in $\Sigma$ if $B$
is chosen sufficiently large in the definition of $\Sigma$.
\begin{lemma}
\label{lem:counting0}
For any $0 \leq \varsigma\leq \alpha_0$ 
and each $W \in \Sigma$
\[
\sum_{W_i \in\W_n} |W_i|^\varsigma ||DT^n|^{-1} J_WT^n|_{\C^0(W_i)}
\leq C \sum_{k=1}^n \delta^{\varsigma-1}\rho^{n-k}\int_{W}|DT^{-k}| +C|W|^{\varsigma}\rho^{n} 
\]
where $J_WT^n$ denotes the Jacobian of $T^n$ along the leaf $T^{-n}W$.
\end{lemma}
\begin{proof}
For each $1 \leq k \leq n$, denote by $W_i^k$ the elements of $\W_k$.
Let $A_k = \{i : |W_i^k|<\delta \}$
and $B_k = \{ i:|W_i^k| \geq \delta \}$ denote the short and long pieces in $\W_k$
respectively.  We regard $\{ W_i^k \}_{i,k}$ as a tree with $W$ as its
root and $\W_k$ as the $k^{\mbox{\scriptsize th}}$ level.

At level $n$, we collect the short pieces into groups as follows.
Consider a piece $W^n_{i_0}\in\W_n$, not necessarily short.  Let $W^k_j$ be the most 
recent long ``ancestor''
of $W^n_{i_0}$, i.e., $k = \max \{0 \leq m \leq n : T^{n-m}(W_{i_0}^n)
\subset W_j^m \; \mbox{and} \; j \in B_m \}$. 
If no such ancestor exists, set $k=0$ and $W^k_j = W$.  
Note that if $W^n_{i_0}$ is long, then $W^k_j = W^n_{i_0}$.
Let
\[
J_n(W^k_j) = \{ i : T^{n-k}(W_i^n) \subset W_j^k \; \mbox{and} \; 
|T^\ell(W_i^n)| < \delta \; \mbox{for} \;  0 \leq \ell \leq n-k-1 \}
\]
be the set of indices corresponding to the short pieces which have the same most
recent long ancestor as $W^n_{i_0}$, or the set $\{W^n_{i_0}\}$ if the
piece is long.  Since for any $i \in J_n(W^k_j) $, $|T^\ell(W_i^n)| < \delta$
for all $0 \leq \ell \leq n-k-1$, we may estimate $\#J_n(W^k_j) \leq L^{n-k}$ using 
Remark~\ref{rem:iterate p1} and Convention~\ref{convention: n_0=1}.
So using the distortion bounds given by equations 
\eqref{eq:distortion} and \eqref{eq:angle dist}, we estimate
\begin{equation}
\label{eq:short pieces}
\begin{split}
\sum_{i \in J_n(W^k_j)}&  |W_i^n|^\varsigma ||DT^n|^{-1} J_WT^n|_{\C^0(W_i^n)} \\ 
  & \leq  C \sum_{i \in J_n(W^k_j)} |T^{n-k}W_i^n|^\varsigma
|(J_WT^{n-k})^{1-\varsigma}  
           |DT^{n-k}|^{-1}|_{\C^0(W_i^n)} ||DT^k|^{-1} J_WT^k|_{C^0(W_j^k)}  \\
  & \leq C ||DT^k|^{-1} J_WT^k|_{C^0(W_j^k)} |W_j^k|^\varsigma 
           (L^{1-\varsigma} \lambda^{-1}\mu^{-\varsigma})^{n-k}\\
&  \leq C ||DT^k|^{-1} J_WT^k|_{C^0(W_j^k)} |W_j^k|^\varsigma  \rho^{n-k} 
\end{split}
\end{equation}
where for the first inequality we have estimated $|W_i^n||J_WT^{n-k}|_{\C^0(W_i^n)} \leq C |T^{n-k}W_i^n|$,
and for the second we have used the H\"{o}lder inequality.
Grouping all $i\in A_n$ in this way, we are left with estimates over long pieces only, so that
using \eqref{eq:short pieces},
\begin{equation}
\label{eq:long remaining}
\begin{split}
\sum_i |W_i^n|^\varsigma ||DT^n|^{-1} J_WT^n|_{\C^0(W_i^n)} 
&=\sum_{k=0}^n\sum_{j\in B_k}\; \sum_{\ i\in J_n(W^k_j)}|W_i^n|^\varsigma
||DT^n|^{-1} J_WT^n|_{\C^0(W_i^n)}\\ 
&\leq C \sum_{k=0}^n \sum_{j\in B_k} |W_j^k|^\varsigma ||DT^k|^{-1}J_WT^k|_{\C^0(W_j^k)} 
 \rho^{n-k} .  
\end{split}
\end{equation}
For each $k \geq 1$, we have $|W_j^k| \geq \delta$ and $T^kW^k_{j_1} \cap T^kW^k_{j_2} = \emptyset$ if
$j_1 \neq j_2$.  So we may sum over $j$, again using \eqref{eq:distortion},
\begin{equation}
\label{eq:switch to integral}
\begin{split}
\sum_{j\in B_k} |W_j^k|^\varsigma ||DT^k|^{-1} J_WT^k|_{\C^0(W_j^k)}
    &\leq C \sum_{j\in B_k}|W_j^k|^{\varsigma-1}
\int_{W_j^k} |DT^k|^{-1} J_WT^k dm \\
&\leq C\delta^{\varsigma-1}\int_{W}|DT^{-k}|\, dm \;. 
\end{split}
\end{equation}
Putting \eqref{eq:switch to integral} together with \eqref{eq:long remaining}, we conclude that
\begin{equation}
\label{eq:summable}
\sum_i |W_i^n|^\varsigma ||DT^n|^{-1} J_WT^n|_{\C^0(W_i^n)} \; \leq \;
   C \sum_{k=1}^n \delta^{\varsigma-1} \rho^{n-k} \int_{W}|DT^{-k}| \, dm  +C|W|^{\varsigma}\rho^{n}
\end{equation}
which proves the lemma.
\end{proof}
As an immediate corollary of the above lemma we have
\begin{lemma}
\label{lem:counting}
For any $0 \leq \varsigma\leq \alpha$ 
and each $W \in \Sigma$
\[
\sum_{W_i \in\W_n} |W_i|^\varsigma ||DT^n|^{-1} J_WT^n|_{\C^0(W_i)}
\leq C D_n \delta^{\varsigma-\alpha}|W|^\alpha +C|W|^{\varsigma}\rho^{n} .
\]
\end{lemma}

Next, we have a fundamental lemma that will allow us to establish a connection
between our Banach spaces and the standard spaces of distributions.
\begin{lemma}
\label{lem:distribution}
For each $h \in \C^1(\M)$, $n\geq 0$, and 
$\vf\in\C^1_{\Si^+_n}$ we have
\[
|h(\vf)| = | \int_{\M} h \vf \, dV | \leq C_\delta |h|_w(|\vf|_\infty + |D^s\vf|_\infty)
\]
where $D^s$ denotes the derivative along the stable direction and $dV$ is the normalized volume
element on $\M$.
\end{lemma}
\begin{proof}
Choose
$\vf\in\C^1_{\Si_{n_0}^+}$ for some $n_0\in\N$, so that $\vf\in\C^1(\overline K)$
for each $K\in\mathcal K_{n_0}$.  
First partition each $\M_i^+$ into finitely many approximate boxes $B_{\ell}$ whose boundary
curves are elements of $\Sigma$ and $\F^u$, as well as the boundary curves of $\M_i^+$ where
necessary. 
The $B_{\ell}$ can be constructed so that each $B_{\ell}$ is foliated by curves  
$W \in \Sigma$ and $\mbox{diam}(B_{\ell}) \leq 2\delta$.
On each $B_{\ell}$, choose a smooth partition $\{ W_{\ell}(\xi) \}$ of $B_{\ell}$
made up of elements of $\Sigma$ which completely cross $B_{\ell}$ in
the approximate stable direction.   
Here $\xi \in E_{\ell}$ is a parameter which indexes
the elements of the foliation $\{ W_{\ell}(\xi) \}$.  

In order to integrate along the curves $W_\ell(\xi)$, we decompose the measure $dV$ on $B_\ell$
into $dV = \nu(d\xi) dm_{\xi, \ell}$ where $m_{\xi,\ell}$ is the
conditional measure on the fiber $W_{\ell}(\xi)$ and $\nu$ is an appropriate measure on 
$\cup_\ell E_\ell$.  We normalize the measures so that
$m_{\xi,\ell}(W_{\ell}(\xi))=|W_{\ell}(\xi)|$; thus, since the
foliation is smooth, $dm_{\xi,\ell}=\rho_{\xi,\ell}dm$ where $m$ is
the arc-length measure on $W_{\ell}(\xi)$ and
$\rho_{\xi,\ell},\rho_{\xi,\ell}^{-1}\in\C^1(W_{\ell}(\xi))$. Note that
$\nu(E_\ell)<\infty$. 

Taking $n\geq n_0$, we estimate
\begin{equation}
\label{eq:decompose}
\begin{split}
\int_{\M} h \vf \; dV & =  \sum_{\ell} \int_{B_{\ell}} \Lp^n h \; \vf \circ T^{-n} \; dV \; 
         = \; \sum_{\ell} \int_{E_{\ell}} \nu(d\xi) \int_{W_{\ell}(\xi)}
\Lp^n h \; \vf \circ T^{-n} \; dm_{\xi,\ell}             \\
      & = \sum_\ell \int_{E_\ell} \nu(d\xi) \int_{W_\ell(\xi)} 
h\circ T^{-n} |DT^n(T^{-n})|^{-1} \vf \circ T^{-n} \rho_{\xi,\ell} \, dm 
\end{split}
\end{equation}
where we have used the definition of $\Lp^n h$ in the second line.  We estimate the integral by
changing variables on one $W_\ell(\xi)$ at a time.
\begin{equation}
\label{eq:one leaf}
\begin{split}
\int_{W_\ell(\xi)} &
\frac{h\circ T^{-n}}{ |DT^n(T^{-n})|}\; \vf \circ T^{-n} \rho_{\xi,\ell} \, dm 
\; = \; \sum_i \int_{W^n_{\ell,i}(\xi)}
h \vf |DT^n|^{-1} J_WT^n \rho_{\xi,\ell}\circ T^{n}\; dm    \\ 
& \leq \; C |h|_w \sum_i
         |\vf|_{\C^1(W^n_{\ell,i}(\xi))} 
         | |DT^n|^{-1} J_WT^n|_{\C^1(W^n_{\ell,i}(\xi))}
         |\rho_{\xi,\ell} \circ T^n|_{\C^1(W^n_{\ell,i}(\xi))}
\end{split}
\end{equation}
where $W^n_{\ell,i}(\xi)$ are the smooth components of
$T^{-n}W_{\ell}(\xi)$ as defined earlier.  Note that 
$|\rho_{\xi,\ell} \circ T^n|_{\C^1(W^n_{\ell,i}(\xi))} \leq C |\rho_{\xi,\ell}|_{\C^1(W_\ell(\xi))} \leq C$ for some
$C$ independent of $\xi$ and $\ell$ using the estimate of equation~\eqref{eq:C1 C0} in Section~\ref{weak norm}.
Also, the distortion bounds of Appendix A imply that
$| |DT^n|^{-1} J_WT^n|_{\C^1(W^n_{\ell,i}(\xi))} \leq | |DT^n|^{-1} J_WT^n|_{\C^0(W^n_{\ell,i}(\xi))}$.

As $n$ increases, elements of $T^{-n}\Sigma$ become more closely aligned with the 
stable direction.  So we may choose an $n_1$, depending on $\vf$, but
not on $\ell$ or $\xi$, such that for $n \geq n_1$ and each $i$,
$|\vf|_{\C^1(W^n_{\ell,i}(\xi))} \leq 2 (|\vf|_\infty +
|D^s\vf|_\infty)$.
For $n \geq n_1+n_0$, using \eqref{eq:one leaf}, we estimate 
\begin{equation}
\label{eq:one leaf estimate}
\int_{W_{\ell}(\xi)} \Lp^n h \, \vf\circ T^{-n} \; dm_{\xi,\ell} 
    \leq \; C | h |_w (|\vf|_\infty + |D^s \vf|_\infty) 
           \sum_i ||DT^n|^{-1} J_WT^n|_{\C^0(W^n_{\ell,i}(\xi))}. 
\end{equation}
To estimate the sum in \eqref{eq:one leaf estimate}, we use Lemma~\ref{lem:counting0}, with $\varsigma = 0$. 
\[
\sum_i ||DT^n|^{-1} J_WT^n|_{\C^0(W_{\ell}^n(\xi))}
   \leq C \sum_{k=1}^n \delta^{-1}\rho^{n-k}\int_{W}|DT^{-k}|dm +C\rho^{n} 
\]

This, together with \eqref{eq:one leaf estimate}, allows us to estimate \eqref{eq:decompose}.
\begin{eqnarray*}
\int_{\M} h \vf \; dV & \leq & C | h |_w (|\vf|_\infty + |D^s \vf|_\infty)     \\       
      & & \; \; \cdot \left( \sum_{\ell} \int_{E_{\ell}} \nu(d\xi)  \; \rho^n    
            + \delta^{-1} \sum_{k=1}^n \int_{\M} |DT^{-k}| dV 
            \; \rho^{n-k} \right).
\end{eqnarray*}

 Since the integral 
$\int_{\M} |DT^{-k}| dV = 1$ for each $k$, the sum
over $k \geq 1$ is bounded independently of $n$.  This proves the lemma.   
\end{proof}


\subsection{Embeddings and Compactness}
\label{sub:embedding}

Notice that, by definition,
$|\cdot|_w\leq C\|\cdot\|_s$. This means that there exists a natural
embedding of $\B$ into $\B_w$. In addition, if $h\in\B$ and $|h|_w=0$,
it is immediate from the definitions \eqref{eq:weak},
\eqref{eq:s-stable} and \eqref{eq:s-unstable} that $\|h\|=0$, i.e.\ that
the embedding is injective.  Accordingly, we will consider
$\B$ as a subset of $\B_w$ in what follows.
\begin{rem}
\label{rem:inclusion}
Lemma \ref{lem:distribution} implies that, for each $h\in\B_w$ and 
$\vf\in\C^1_{\Si_n^+}$,
$|h(\vf)|\leq C|h|_w|\vf|_{\C^1_{\Si_n^+}}$, that is,
$\C^1\hookrightarrow\B\hookrightarrow\B_w\hookrightarrow
(\C^1_{\Si_n^+})'$. 
In fact, the inclusions 
are injective: 
if $h_1,h_2$ coincide as 
elements of $(\C^1_{\Si_n^+})'$ and they both belong to any of the spaces $\C^1$, $\B$, or $\B_w$, then
they coincide as elements of those spaces as well.
This can be proven as
in \cite[Proposition 4.1]{liverani gouezel}.
\end{rem}

We can finally state the last result of this section.
\begin{lemma}
\label{lem:compact}
The unit ball of $\B$ is compactly
embedded in $\B_w$. 
\end{lemma}

To prove the above fact it is convenient to remark the following
obvious result.
\begin{lemma}
\label{lem:compact leaf}
For any fixed $W \in \Sigma$, the unit ball of $| \cdot |_{\C^1(W)}$ 
is compactly embedded in $| \cdot |_{W, \alpha, q}$.
\end{lemma}
\begin{proof}  For fixed $W$, 
$| \cdot |_{W, \alpha, q}$ is equivalent to $| \cdot |_{\C^q(W)}$.  
The lemma follows immediately.
\end{proof}

\begin{proof}[\bf Proof of Lemma \ref{lem:compact}]
Since on each leaf $W \in \Sigma$, $\|\cdot\|_s$ is the dual of $| \cdot |_{W, \alpha, q}$ 
and $|\cdot|_w$ is the dual of 
$| \cdot |_{\C^1(W)}$, Lemma~\ref{lem:compact leaf} implies that the unit ball of 
$\|\cdot\|_s$ is compactly embedded in 
$|\cdot |_w$ on $W$.  It remains to compare the weak norm on different leaves.

Let $0 < \ve \leq \ve_0$ be fixed.  The set of functions $\Xi$ is compact in the $\C^1$-norm
so on each $\M_i^+$, we may choose finitely many leaves $W^i \in \Sigma$ such that $\{ W^i \}$ 
forms an $\ve$-covering
of $\Sigma|_{\M_i}$ in the distance $d_\Sigma$. 
Since any ball of finite radius in the $\C^1$-norm is compactly
embedded in $\C^q$, we may choose finitely many functions 
$\overline{\vf}_j \in \C^1(I_{r_0})$ such that $\{ \overline{\vf}_j \}$ forms an 
$\ve$-covering in the $\C^q(I_{r_0})$-norm of the ball of
radius $(1+\kappa)^2$ in $\C^1(I_{r_0})$.

Now let $h \in \C^1(\M)$, $W \in \Sigma$, and $\vf \in \C^1(W)$ with
$|\vf|_{\C^1(W)} \leq 1$. 
Let $F$ denote the function associated with $W$ and as usual, let
$G_F$ be the lift of the graph of $F$ to $\M$. 
Let $\overline{\vf} = \vf \circ G_F$ be the push down of $\vf$ to $I_r$.  Note that 
$|\overline{\vf}|_{\C^1(I_r)} \leq (1+\kappa)^2$.

Choose $W^i$ such that $d_\Sigma(W,W^i) \leq \ve$ and $\overline{\vf}_j$ such that 
$|\overline{\vf} - \overline{\vf}_j|_{\C^q(I_r)} \leq \ve$.
Let $F^i$ and $G_{F^i}$ denote the usual functions associated with the
leaf $W^i$ and define $\vf_j = \overline{\vf}_j \circ G_{F^i}^{-1}$. 
Note that $|\vf_j|_{\C^1(W^i)} \leq (1+\kappa)^3$.
Then normalizing $\vf$ and $\vf_j$ by $(1+\kappa)^3$, we get
\[
\left| \int_W h \vf \; dm - \int_{W^i} h \vf_j \; dm \right|
       \; \leq \; \ve^\beta \|h\|_u (1+\kappa)^3  \;
       \leq \; \ve^\beta (1+\kappa)^3 b^{-1} \|h\| . 
\]
We have proved that for each $0 < \ve \leq \ve_0$, there exist finitely many
bounded linear functionals $\ell_{i,j}$, $\ell_{i,j}(h) = \int_{W^i} h \vf_j dm$, such that
\[
|h|_w \leq \sup_{i,j} \ell_{i,j}(h) + \ve^\beta C_b \|h\|
\]
which implies the desired compactness.
\end{proof}


\section{Lasota-Yorke Estimates}
\label{inequalities}

In this section we prove Proposition \ref{prop:lasotayorke}.


\subsection{Estimating the Weak Norm}
\label{weak norm}

For $h \in \C^1(\M)$, $W \in \Sigma$ and $\vf \in \C^1(W)$ such that $|\vf|_{\C^1(W)} \leq 1$, we have
\begin{equation}
\label{eq:start}
\int_W\Lp^nh \; \vf \; dm =\int_{T^{-n}W}h\frac{J_WT^n}{|DT^n|}\vf\circ T^n dm
=\sum_{W_i \in\W_n}\int_{W_i}h\frac{J_WT^n}{|DT^n|}\vf\circ T^n dm
\end{equation}
where as before $J_WT^n$ denotes the Jacobian of $T^n$ along the leaf $T^{-n}W$.

Using the definition of the weak norm on each $W_i$, we estimate \eqref{eq:start} by 
\begin{equation}
\label{eq:weak estimate}
\int_W\Lp^nh \; \vf\; dm \; \leq \; \sum_{W_i \in\W_n} |h|_w
| |DT^n|^{-1}J_WT^n|_{\C^1(W_i)} |\vf\circ T^n|_{\C^1(W_i)} .
\end{equation}
The disortion bounds given by equation \eqref{eq:distortion} imply that
\[ 
| |DT^n|^{-1}J_WT^n|_{\C^1(W_i)} \leq C | |DT^n|^{-1}J_WT^n|_{\C^0(W_i)} .
\]
Also notice that 
\begin{equation}
\label{eq:C1 C0}
\frac{|\vf (T^nx) - \vf (T^ny)|}{d_s(T^nx,T^ny)}\cdot \frac{d_s(T^nx,T^ny)}{d_s(x,y)} \leq C |\vf|_{\C^1(W)} |J_WT^n|_{\C^0(W_i)} 
    \leq C\lambda^{-n} |\vf|_{\C^1(W)}
\end{equation}
for any $x,y \in W_i$,
so that $|\vf \circ T^n|_{\C^1(W_i)} \leq C |\vf|_{\C^1(W)}$.
Using these estimates in equation \eqref{eq:weak estimate}, we obtain
\[
\int_W \Lp^n h \; \vf \; dm \; \leq \; C |h|_w |\vf|_{\C^1(W)} \sum_{W_i \in\W_n}
||DT^n|^{-1} J_WT^n|_{\C^0(W_i)}  . 
\]
The above formula, together with  Lemma \ref{lem:counting} used in the case
$\varsigma=0$, yields the inequality,
\[
\int_W \Lp^nh \; \vf \; dm \leq C |h|_w (D_n + \rho^n) |\vf|_{\C^1(W)}.
\]
Taking the supremum over all $W \in \Sigma$ and $\vf \in \C^1(W)$ with $|\vf|_{\C^1(W)} \leq 1$ 
yields the required estimate \eqref{eq:weak norm}.


\subsection{Estimating the Strong Stable Norm}
\label{stable norm}

Using equation \eqref{eq:start}, we write for each $W\in\Sigma$ and
$\vf\in\C^1(W,\mathbb{C})$ such that $|\vf|_{W,\alpha,q}\leq 1$,
\begin{equation}
\label{eq:stable split}
\int_W\Lp^nh\; \vf \; dm  =  \sum_{i}\left\{\int_{W_i}h\frac{J_WT^n}{|DT^n|}\,\overline{\vf}_i
 +\frac 1{|W_i|}\int_{W_i}\vf\circ T^n  \int_{W_i}h\frac{J_WT^n}{|DT^n|} \right\},
\end{equation}
where $\overline{\vf}_i:=\vf\circ T^n -\frac 1{|W_i|}\int_{W_i}\vf\circ T^n \, dm$.
Let us estimate the above expression.

To estimate the first term of \eqref{eq:stable split},
we first estimate $|\overline{\vf}_i|_{\C^q(W_i)}$. 

Following equation~\eqref{eq:C1 C0}, we write
\begin{equation}
\label{eq:H^q}
\frac{|\vf (T^nx) - \vf (T^ny)|}{d_s(x,y)^q} \leq C |J_WT^n|^q_{\C^1(W_i)} |\vf|_{\C^1(W)}
\end{equation}
for any $x,y \in W_i$.  If $H^q(f)$ represents the H\"{o}lder constant of $f$, then
\eqref{eq:H^q} implies that $H^q(\vf \circ T^n) \leq C |J_WT^n|^q_{\C^0(W_i)} H^q(\vf)$ on approximate stable leaves
due to the contraction of $T^n$.
Also,
\[
\begin{split}
& \left| \vf \circ T^n - \frac{1}{|W_i|} \int_{W_i} \vf \circ T^n \, dm \right|_{\C^0(W_i)} 
    \leq |\sup_{W_i} \vf \circ T^n - \inf_{W_i} \vf \circ T^n|  \\
    & \qquad \qquad \leq H^q(\vf \circ T^n) |W_i|^q \leq C H^q(\vf) |J_WT^n|^q_{\C^0(W_i)}  .
\end{split}
\]
This estimate together with \eqref{eq:H^q} and the fact that $|\vf|_{W,\alpha,q} \leq 1$, implies
\begin{equation}
\label{eq:C^q small}
|\overline{\vf}_i|_{\C^q(W_i)} \leq C |J_WT^n|_{\C^0(W_i)}^q |\vf|_{\C^q(W)} \leq C| J_WT^n|_{\C^0(W_i)}^q |W|^{-\alpha} .
\end{equation}

Applying \eqref{eq:C^q small} and the definition of the strong stable norm to the first term of 
\eqref{eq:stable split} yields,
\begin{equation}
\label{eq:first stable}
\begin{split}
\sum_{i} \int_{W_i}h & \frac{J_WT^n}{|DT^n|}  \, \overline{\vf}_i \; dm  \leq 
C\sum_{i} \|h\|_s |W_i|^\alpha ||DT^n|^{-1} J_WT^n|_{\C^0(W_i)} |J_WT^n|_{\C^0(W_i)}^q |W|^{-\alpha} \\
& \leq \; C \|h\|_s |W|^{-\alpha} \mu_+^{qn} \sum_i |W_i|^\alpha ||DT^n|^{-1} J_WT^n|_{\C^0(W_i)}
\; \leq \; C D_n\|h\|_s\mu_+^{qn}.          
\end{split}
\end{equation}
where in the second line we have used Lemma~\ref{lem:counting} with $\varsigma = \alpha$.

For the second term of \eqref{eq:stable split}, we use the fact that $|\vf|_{\infty} \leq |W|^{-\alpha}$ 
to estimate $\frac 1{|W_i|}\int_{W_i}\vf \circ T^n dm \leq |W|^{-\alpha}$. 
Recall the notation used in the proof of Lemma
\ref{lem:counting0}.  Grouping the pieces $W_i = W^n_i$ according to most recent long ancestors, we have
\[
\begin{split}
\sum_{i}  |W|^{-\alpha} \int_{W_i}h|DT^n|^{-1}J_WT^n \, dm
= & \sum_{k=1}^n\sum_{j\in B_k}\sum_{\ i\in J_n(W^k_j)} 
     |W|^{-\alpha} \int_{W_i}h|DT^n|^{-1}J_WT^n \, dm \\
&  + \sum_{\ i\in J_n(W^0_j)} 
     |W|^{-\alpha} \int_{W_i}h|DT^n|^{-1}J_WT^n \, dm
\end{split}
\]
where we have split up the terms involving $k=0$ and $k \geq 1$.  
We estimate the terms with $k \geq 1$
by the weak norm and the terms with $k=0$ by the strong stable norm,
\[
\begin{split}
\sum_{i}  |W|^{-\alpha}\left| \int_{W_i}h|DT^n|^{-1}J_WT^n \, dm\right|
\leq &
C\sum_{k=1}^n\sum_{j\in B_k}\sum_{\ i\in
J_n(W^k_j)}|W|^{-\alpha}|h|_w ||DT^n|^{-1} J_WT^n|_{\C^0(W_i)}\\
&  +
C \sum_{\ i\in J_n(W^0_j)} |W|^{-\alpha} \|h\|_s |W_i|^{\alpha} ||DT^n|^{-1} J_WT^n|_{\C^0(W_i)} .
\end{split}
\]
Using equations \eqref{eq:switch to integral} and \eqref{eq:summable} from Lemma \ref{lem:counting0}, 
with $\varsigma = 0$ for the first sum
and $\varsigma = \alpha$ for the second, we conclude

\begin{equation}
\label{eq:second stable}
\sum_{i}  \frac 1{|W_i|}\left|\int_{W_i}\vf \circ T^n \, dm \int_{W_i}h|DT^n|^{-1}J_WT^n \, dm\right|
\; \leq \; CD_n \delta^{-\alpha}|h|_w + C\|h\|_s\rho^{n} \\
\end{equation}

Putting together \eqref{eq:first stable} and \eqref{eq:second stable} proves \eqref{eq:stable norm},
\[
\|\Lp^n h\|_s \leq C\left( D_n\mu_+^{q n}+\rho^n\right)\|h\|_s 
         + C_\delta  D_n |h|_w .
\]


\subsection{Estimating the Strong Unstable Norm}
\label{unstable norm}

Consider two admissible leaves $W^i\in\Sigma$, $d_\Sigma(W^1,W^2) \leq \ve$.  They can be
partitioned into ``matched'' pieces $U^i_j$ and
``unmatched'' pieces $V^i_j$. To do so consider the connected pieces of
$W^i\setminus\Si^-_n$. If one looks at their image under $T^{-n}$ then
one can associate to each point $x \in T^{-n}(W^1\cup W^2)$ a vertical (in the chart) segment
$\gamma_x \in \F^u$, of length at most $C\lambda^{-n}\ve$, such that its image under $T^n$, 
if not cut by a singularity, will be of length $C\ve$ centered at $x$. We can thus
subdivide the connected pieces of $W^i\setminus\Si^-_n$ into subintervals of points for which
$T^n\gamma_x$ intersects the other manifold and subintervals for which this
is not the case. In the latter case, we call the subintervals $V^i_j$ and note that 
either we are at the endpoints of $W^i$ or the vertical segment is cut by a singularity. In
both cases the subintervals $V^i_j$ can be of length at most $C\ve$ and their number
is at most $L_n+2$.\footnote{Without any loss of information (by throwing out at most finitely
many points), we can take each $V^i_j$ to be the image of an open interval. Thus for fixed $i$, the
$V^i_j$ are disjoint.} 
In the remaining
pieces the curves $T^n\gamma_x$ provide a one to one correspondence 
between points in $W^1$ and $W^2$. We can further partition the pieces in
such a way that the lengths of their preimages are between $\delta$ and $2\delta$ and the
partitioning can be made so that the pieces are pairwise matched by
the foliation $\{\gamma_x\}$. We call these matched pieces $U^i_j$.
In this way we write $W^i = (\cup_j U^i_j) \cup (\cup_k V^i_k)$.
Note that the unmatched pieces $V^i_j$ must be short while the matched pieces $U^i_j$
may be long or short.

To be more precise, remember that to exactly describe the leaf
$T^{-n}U^1_j$ we must give $i_j, x_j, r_j,F^1_j$ 
so that $T^{-n}U^1_j=\chi_{i_j}(G(x_j, r_j,F^1_j)(I_{r_j}))$ (see the
end of section \ref{leaves}). Once the leaves $T^{-n}U^1_j$ are
described in such a way we have, by construction, that
$T^{-n}U^2_j$ is of the form $G(x_j,
r_j,F^2_j)(I_{r_j})$ for some appropriate function $F^2_j$ so that the
point $z:=x_j+(t,F^1_j(t))$ is associated with the point
$x_j+(t,F^2_j(t))\in\chi_{i_j}^{-1}(T^{-n}U^2_j)$ by the vertical
segment $\chi_{i_j}^{-1}(\gamma_{\chi_{i_j}(z)})=\{(0,s)\}_{s\in\R}$.

Given $\vf_i$ on $W^i$ with $|\vf_i|_{\C^1(W^i)} \leq 1$ and $d_q(\vf_1, \vf_2) \leq \ve$, 
with the above construction we must estimate
\begin{equation}
\label{eq:unstable split}
\begin{split}
& \left|\int_{W^1} \Lp^nh \, \vf_1 \, dm - \int_{W^2} \Lp^nh \, \vf_2 \, dm\right| 
  \; \leq \; \sum_{i,j} \left|\int_{T^{-n}V^i_j} h |DT^n|^{-1}J_{W^i}T^n \vf_i\circ T^n \, dm \right|\\
  & + \sum_j \left| \int_{T^{-n}U^1_j} h |DT^n|^{-1}J_{W^1}T^n \vf_1\circ T^n \, dm 
    - \int_{T^{-n}U^2_j} h |DT^n|^{-1}J_{W^2}T^n \vf_2\circ T^n \, dm \right|
\end{split}
\end{equation}

We do the estimate over the unmatched short pieces $V^i_j$ first.
Using the strong stable norm, we can compute,
\begin{equation}
\label{eq:first unstable}
\begin{split}
\sum_{i,j} & \left|\int_{T^{-n}V^i_j}h|DT^n|^{-1}J_WT^n\vf_i\circ T^n \, dm\right| \leq 
\sum_{i,j} \|h\|_s |T^{-n}V^i_j|^\alpha ||DT^n|^{-1}J_WT^n|_{\C^q} |\vf_i|_{\C^q}\\
& \leq C \| h \|_s \sum_{i,j} |V^i_j|^\alpha ||DT^n|^{-1} (J_WT^n)^{1-\alpha}|_{\C^0} 
\leq C \ve^\alpha \|h\|_s L_n \lambda^{-n} \mu^{-\alpha n} .
\end{split}
\end{equation}

Next, we must estimate
\[
\sum_j\left|\int_{T^{-n}U^1_j}h|DT^n|^{-1}J_{W^1}T^n \, \vf_1\circ T^n \, dm -
\int_{T^{-n}U^2_j}h|DT^n|^{-1}J_{W^2}T^n \, \vf_2\circ T^n \, dm \right|   .
\]
First, recall that for each $F\in\Xi$, $G_F(t) = \chi(x_F + (t,F(t))$ for $t \in I_r$.
We define the map $\Psi :U^2_j\to U^1_j$ by $\Psi:=T^n \circ
G_{F^1_j}\circ G_{F^2_j}^{-1} \circ T^{-n}$  and the function
\[
\tilde\vf:=\left[\vf_1\cdot (|DT^n|^{-1}J_{W^1}T^n)\circ T^{-n}\right]\circ\Psi \cdot 
[(|DT^n|^{-1}J_{W^2}T^n)\circ T^{-n}]^{-1}.
\]
$\tilde \vf$ is well-defined on $U^2_j$ and 
$\left[\tilde\vf\circ T^n |DT^n|^{-1}J_{W^2}T^n\right]
\circ G_{F^2_j}=\left[\vf_1\circ T^n|DT^n|^{-1}J_{W^1}T^n\right]\circ G_{F^1_j}$.

We can then write
\begin{equation}
\label{eq:stepone}
\begin{split}
&\sum_j\left|\int_{T^{-n}U^1_j}h|DT^n|^{-1}J_{W^1}T^n\vf_1\circ T^n-
\int_{T^{-n}U^2_j}h|DT^n|^{-1}J_{W^2}T^n\vf_2\circ T^n\right|\\
&\leq \sum_j\left|\int_{T^{-n}U^1_j}h|DT^n|^{-1}J_{W^1}T^n\vf_1\circ T^n-
\int_{T^{-n}U^2_j}h|DT^n|^{-1}J_{W^2}T^n\tilde\vf\circ T^n\right|\\
&\quad+\sum_j\left|\int_{T^{-n}U^2_j}h|DT^n|^{-1}J_{W^2}T^n(\tilde\vf-\vf_2)\circ
T^n\right| .
\end{split}
\end{equation}

We estimate the first sum in equation~\eqref{eq:stepone} using the strong unstable norm.

The distortion bounds given by \eqref{eq:distortion} and the estimate of \eqref{eq:C1 C0} imply that
\begin{equation}
\label{eq:c1-unst 1}
|\,|DT^n|^{-1}J_{W^1}T^n\cdot \vf_1 \circ
T^n|_{\C^1(T^{-n}U^1_j)}\leq 
C|\vf_1|_{\C^1(U^1_j)}|\,|DT^n|^{-1}J_{W^1}T^n|_{\C^0(T^{-n}U^1_j)}  .
\end{equation}
Also, by the definition of $\tilde\vf$,
\begin{equation}
\label{eq:c1-unst 2}
\begin{split}
& \left|\,\frac{J_{W^2}T^n}{|DT^n|}\, \tilde\vf \circ T^n\right|_{\C^1(T^{-n}U^2_j)} 
  = \left |\left[\vf_1 \circ T^n  \cdot
    (|DT^n|^{-1}J_{W^1}T^n)\right]\circ G_{F^1_j} \circ G^{-1}_{F^2_j} \right|_{\C^1(T^{-n}U^2_j)}\\ 
&\quad\quad\leq C|\vf_1|_{\C^1(U^1_j)}|\,|DT^n|^{-1}J_{W^1}T^n|_{\C^0(T^{-n}U^1_j)}   .
\end{split}
\end{equation}
By the definition of $\tilde\vf$ and $d_q$,
\[
\begin{split}
d_q(|DT^n|^{-1} & J_{W^2}T^n\tilde\vf\circ T^n,|DT^n|^{-1}J_{W^1}T^n\vf_1\circ T^n) \\
  & = \left| \left[ |DT^n|^{-1} J_{W^2}T^n\tilde\vf\circ T^n \right] \circ G_{F^2_j}
    - \left[  |DT^n|^{-1}J_{W^1}T^n\vf_1\circ T^n \right] \circ G_{F^1_j} \right| \; = \; 0 .
\end{split}
\]
In addition, the uniform hyperbolicity of $T$ implies that
\[
d_\Sigma(T^{-n}U^1_j,T^{-n}U^2_j)\leq C\lambda^{-n}\ve =: \ve_1 .
\]  
This follows
from the usual graph transform argument which is standard to hyperbolic theory.

We first renormalize the test functions by 
$R_j = C|\vf_1|_{\C^1}\,||DT^n|^{-1}J_{W^1}T^n|_{\C^0(T^{-n}U^1_j)}$.
Then for each $j$, we apply the definition of the strong unstable norm with $\ve_1$ in place of $\ve$.
Thus,

\begin{equation}
\label{eq:second unstable}
\begin{split}
\sum_{j} &\left|\int_{T^{-n}U^1_j}h  |DT^n|^{-1} J_{W^1}T^n \, \vf_1\circ T^n -
\int_{T^{-n}U^2_j} h  |DT^n|^{-1} J_{W^2}T^n \, \tilde\vf\circ T^n\right|  \\
& \leq C \ve^\beta_1\sum_{j}  |\vf_1|_{\C^1}\,||DT^n|^{-1}J_{W^1}T^n|_{\C^0(T^{-n}U^1_j)}  \|h\|_u  \\
&   \leq C \|h\|_u \lambda^{-n\beta} \ve^\beta (D_n+\rho^n),
\end{split}
\end{equation}
where we have used Lemma~\ref{lem:counting} in the last line with $\varsigma = 0$.

It remains to estimate the second sum in \eqref{eq:stepone} using the strong stable norm.
We need the following fact.
\begin{lemma}
\label{lem:distortion}
For each $j$, we have that
\[
|(|DT^n|^{-1} J_{W^1}T^n)\circ G_{F^1_j}
-(|DT^n|^{-1} J_{W^2}T^n)\circ G_{F^2_j}|_{\C^q}\leq C | |DT^n|^{-1}
J_{W^1}T^n|_{C^0(T^{-n}U^1_j)} \ve^{1-q}. 
\]
\end{lemma}
\begin{proof}
Throughout the proof, for ease of notation 
we write $J_i^n$ for $|DT^n|^{-1} J_{W^i}T^n$.

For any $t \in I_{r^2_j}$, $x=G_{F^1_j}(t)$ and $y=G_{F_j^2}(t)$ lie
on a common element $\gamma \in \F^u$. 
Thus $T^n(x)$ and $T^n(y)$ also lie on the element $T^n\gamma \in
\F^u$ which intersects $W^1$ and $W^2$ and 
has length at most $C\varepsilon$.
By \eqref{eq:distortion},
\begin{equation}
\label{eq:J C0}
|J_1^n(x) - J_2^n(y)| 
\leq C|J_1^n|_{\C^0(T^{-n}U^1_j)} d_u (T^nx,T^ny) \leq C \ve |J_1^n|_{\C^0(T^{-n}U^1_j)}  .
\end{equation}
Using this estimate and the fact that $J_i^n \circ G_{F^i_j}$ is $\C^1$, we write
\begin{equation}
\label{eq:first holder}
\frac{|(J_1^n \circ G_{F^1_j}(s) - J_2^n \circ G_{F^2_j}(s))
- (J_1^n \circ G_{F^1_j}(t) - J_2^n \circ G_{F^2_j}(t)) |}{|s-t|^q}
\leq \frac{2 C \ve |J_1^n|_{\C^0(T^{-n}U^1_j)} }{|s-t|^q} .
\end{equation}
Also, 
\begin{equation}
\label{eq:second holder}
\frac{|(J_1^n \circ G_{F^1_j}(s) - J_1^n \circ G_{F^1_j}(t))
- (J_2^n \circ G_{F^2_j}(s) - J_2^n \circ G_{F^2_j}(t)) |}{|s-t|^q}
\leq 2 C|J_1^n|_{\C^0(T^{-n}U^1_j)} |s-t|^{1-q}.
\end{equation}
Putting \eqref{eq:first holder} and \eqref{eq:second holder} together implies
that the H\"{o}lder constant of $J_1^n \circ G_{F^1_j}- J_2^n \circ G_{F^2_j}$
is bounded by 
\[
H^q(J_1^n \circ G_{F^1_j}- J_2^n \circ G_{F^2_j}) 
\leq C |J_1^n|_{\C^0(T^{-n}U^1_j)} \min \{ \ve |s-t|^{-q}, |s-t|^{1-q} \} .
\]
This expression is maximized when $\ve |s-t|^{-q} = |s-t|^{1-q}$, i.e., when
$\ve = |s-t|$.  Thus $H^q(J_1^n \circ G_{F^1_j}- J_2^n \circ G_{F^2_j}) 
\leq C |J_1^n|_{\C^0(T^{-n}U^1_j)} \ve^{1-q}$, which,  
together with \eqref{eq:J C0}, concludes the proof of the lemma.
\end{proof}

Using the strong stable norm, we estimate the second sum in \eqref{eq:stepone} by
\begin{equation}
\label{eq:unstable strong}
\begin{split}
\sum_j & \left|\int_{T^{-n}U^2_j}h|DT^n|^{-1}J_{W^2}T^n(\tilde\vf-\vf_2)\circ T^n \right|   \\
& \leq \; C \|h\|_s \sum_j |T^{-n}U^2_j|^\alpha 
       \left||DT^n|^{-1} J_{W^2}T^n(\tilde\vf-\vf_2)\circ T^n\right|_{\C^q(T^{-n}U^2_j)}
\end{split}
\end{equation}
In order to estimate the $\C^q$-norm of the function in \eqref{eq:unstable strong}, we split
it up into two differences, 
\begin{equation}
\label{eq:diff}
\begin{split}
|& |DT^n|^{-1}J_{W^2}T^n\cdot(\tilde\vf - \vf_2)\circ T^n|_{\C^q(T^{-n}U^2_j)}\\
\leq& \; C\left |\left[ (|DT^n|^{-1}J_{W^1}T^n)\cdot\vf_1 \circ
T^n\right]\circ G_{F^1_j} - \left[(|DT^n|^{-1}J_{W^2}T^n)\cdot\vf_2
\circ T^n\right] \circ G_{F^2_j}\right|_{\C^q(I_{r_j})}\\ 
\leq& \; C\left | (|DT^n|^{-1}J_{W^1}T^n)\circ G_{F^1_j} \left[(\vf_1
\circ T^n\circ G_{F^1_j} -\vf_2 \circ T^n\circ
G_{F^2_j}\right]\right|_{\C^q(I_{r_j})}\\ 
&+ C\left|\left[(|DT^n|^{-1}J_{W^1}T^n) \circ
G_{F^1_j}-(|DT^n|^{-1}J_{W^2}T^n) \circ G_{F^2_j}\right]\vf_2\circ
T^n\circ G_{F^2_j}\right|_{\C^q(I_{r_j})}\\ 
\leq& \; C |\, |DT^n|^{-1}J_{W^1}T^n|_{\C^0(T^{-n}U^1_j)} \left|\vf_1 \circ T^n\circ
G_{F^1_j} -\vf_2 \circ T^n\circ G_{F^2_j}\right|_{\C^q(I_{r_j})}\\ 
&+ C\left|(|DT^n|^{-1}J_{W^1}T^n) \circ
G_{F^1_j}-(|DT^n|^{-1}J_{W^2}T^n) \circ
G_{F^2_j}\right|_{\C^q(I_{r_j})} 
\end{split}
\end{equation}
Note that the second term can be bounded using Lemma
\ref{lem:distortion}. To bound the first term, let $F^i\in\Xi$ be the
function defining $W^i$. Then setting $\alpha_j:=G_{F^2}^{-1}\circ
T^n\circ G_{F^2_j}$, we have that $|\alpha_j|_{\C^q}\leq C$ and
\begin{equation}
\label{eq:vf-unstable}
\begin{split}
& \left|\vf_1 \circ T^n\circ G_{F^1_j} -\vf_2 \circ T^n\circ G_{F^2_j}\right|_{\C^q(I_{r_j})}  \\
& \qquad \qquad = \left| \vf_1 \circ T^n\circ G_{F^1_j} \circ \alpha_j^{-1} \circ \alpha_j 
    -\vf_2 \circ T^n\circ G_{F^2_j} \circ \alpha_j^{-1} \circ \alpha_j \right|_{\C^q(I_{r_j})} \\
& \qquad \qquad = \left| \vf_1 \circ \Psi \circ G_{F^2} \circ \alpha_j 
    -\vf_2 \circ G_{F^2} \circ \alpha_j \right|_{\C^q(I_{r_j})} \\
& \qquad \qquad \leq C\left|\vf_1 \circ \Psi\circ G_{F^2} 
    -\vf_2 \circ G_{F^2}\right|_{\C^q(I_{r_j})}  \\
& \qquad \qquad \leq C\left|\vf_1 \circ \Psi\circ G_{F^2} 
    -\vf_1 \circ G_{F^1}\right|_{\C^q(I_{r_j})}
    + C \left| \vf_1 \circ G_{F^1} - \vf_2 \circ G_{F^2} \right|_{\C^q(I_{r_j})}  \\
& \qquad \qquad \leq C\left|\vf_1 \circ G_{F^1}\circ G_{F^1}^{-1}\circ \Psi\circ G_{F^2} 
    -\vf_1 \circ G_{F^1}\right|_{\C^q(I_{r_j})}+Cd_q(\vf_1,\vf_2).
\end{split}
\end{equation}
Thus we need the following final estimate.
\begin{lemma}
\label{lem:change leaf}
For a fixed $U_j^2$, let $J \subset I_{r_2}$ be an
interval on which $G_{F^1}^{-1} \circ \Psi \circ G_{F^2}$ is defined. Then
\[
|Id - G_{F^1}^{-1} \circ \Psi \circ G_{F^2}|_{\C^1(J)} \leq C\ve.
\]
\end{lemma}
\begin{proof}  
Recall that $\Psi = T^n \circ G_{F_j^1} \circ G_{F_j^2}^{-1} \circ T^{-n}$.
The function $\phi_j := G_{F_j^1} \circ G_{F_j^2}^{-1}$ maps a 
point $x \in T^{-n}U^2_j$ to a point $y \in T^{-n}U^1_j$ which lies on an
curve $\gamma \in \F^u$ containing both $x$ and $y$.  Thus $\Psi$ maps $T^n(x)$ to $T^n(y)$
and these two points lie on $T^n\gamma \in \F^u$.  By the transversality of the family $\F^u$, 
this implies that $d_u(T^n x,\Psi( T^n x)) \leq C\ve$ where $d_u$ denotes distance along curves in
$\F^u$.  Then
\[
\begin{split}
& |Id - G_{F^1}^{-1} \circ \Psi \circ G_{F^2}|_{\C^0(J)}
\; = \; |G_{F^1}^{-1} \circ G_{F^1} - G_{F^1}^{-1} \circ \Psi \circ G_{F^2}| 
\; \leq \; |G_{F^1}^{-1}|_{\C^1} |G_{F^1} - \Psi \circ G_{F^2}|   \\
& \leq (1+\kappa) (|G_{F^1} - G_{F^2}| + |G_{F^2} - \Psi \circ G_{F^2}|)
\; \leq \; (1+\kappa) (\ve + C\ve).
\end{split}
\]
Closeness in the $\C^1$-norm follows from the fact that all the functions involved are
bounded in $\C^2$-norm,
$|G_{F^1} - G_{F^2}|_{\C^1} \leq \ve$, and
\[
|\pd \Psi - 1| = | \pd (T^n \circ \phi_j \circ T^{-n}) - 1| 
= \left| \frac{J_WT^n(\phi_j \circ T^{-n})}{J_WT^n(T^{-n})} \pd \phi_j -1 \right| 
\leq C \ve
\]
where $\pd$ denotes differentiation along $T^{-n}W^2$ and in the last inequality we have used distortion
estimate \eqref{eq:distortion}. 
\end{proof}

We can now estimate equation~\eqref{eq:vf-unstable} using Lemma~\ref{lem:change leaf},
\[
\begin{split}
& \left|\vf_1 \circ T^n\circ G_{F^1_j} -\vf_2 \circ T^n\circ G_{F^2_j}\right|_{\C^q(I_{r_j})} \\
& \qquad \qquad \leq C \left|\vf_1 \circ G_{F^1}\circ ( G_{F^1}^{-1}\circ \Psi\circ G_{F^2} - Id ) \right|_{\C^q(I_{r_j})}
    + Cd_q(\vf_1,\vf_2) \\
& \qquad \qquad \leq C | \vf_1 \circ G_{F^1}|_{\C^q(I_{r_j})} |Id - G_{F^1}^{-1} \circ \Psi \circ G_{F^2}|_{\C^1(I_{r_j})} + C\ve \\
& \qquad \qquad \leq C\ve^{1-q}+C\ve.
\end{split}
\]
The above, together with equation~\eqref{eq:diff} and Lemma \ref{lem:distortion}, implies
\begin{equation}
\label{eq:small holder}
| \, |DT^n|^{-1}J_{W^2}T^n\cdot(\tilde\vf - \vf_2)\circ
T^n|_{\C^q(T^{-n}U^2_j)}\leq C\ve^{1-q}|\,
|DT^n|^{-1}J_{W^1}T^n|_{\C^0(T^{-n}U^2_j)}. 
\end{equation}

Since $1-q \geq \beta$, we can use \eqref{eq:small holder} and \eqref{eq:unstable strong} 
to estimate the second sum of \eqref{eq:stepone} by 
\begin{equation}
\label{eq:third unstable}
\begin{split}
\sum_j & \left|\int_{T^{-n}U^2_j}h|DT^n|^{-1}J_{W^2}T^n(\tilde\vf-\vf_2)\circ T^n \right|   \\
& \leq \; C \|h\|_s \ve^{1-q} \sum_j |T^{-n}U^2_j|^\alpha 
       | \, |DT^n|^{-1} J_{W^2}T^n|_{\C^0(T^{-n}U^2_j)}  \\
& \leq \; C \|h\|_s D_n  |W^2|^\alpha  \ve^\beta
\end{split}
\end{equation}
where in the last line we have again used Lemma~\ref{lem:counting} with $\varsigma = \alpha$.

Combining the estimates from equations \eqref{eq:first unstable}, \eqref{eq:second unstable},
and \eqref{eq:third unstable}, we obtain
\[
\| \Lp^n h \|_u \leq C \|h\|_u \lambda^{-\beta n}D_n 
+ C\|h\|_s (D_n + L_n\lambda^{-n}\mu^{-\alpha n}).
\]
This completes the proof of \eqref{eq:unstable norm}.


\section{Spectral Picture}
\label{spectrum}

From the Lasota-Yorke estimates \eqref{eq:almost-ly} and the
compactness it follows by the standard Hennion argument (see
\cite{Bal} for details) that the spectral 
radius of $\Lp$ is bounded by $(D_N)^{\frac 1N}$ and the essential
spectral radius by $\tau (D_N)^{\frac 1N}$ where one can take $N$
arbitrary large provided $b$ is chosen sufficiently small. But since norms
with different $b$ are all equivalent, the spectral radii are
insensitive to the choice of $b$.  Accordingly, fixing $b$ small enough
once and for all, we see that the spectral 
radius of $\Lp$ is bounded by $D_*:=\limsup \limits_{n \to
\infty}D_n^{1/n} $ and the essential 
spectral radius is bounded by $\tau D_*$. To proceed we need an
estimate of $D_*$.


\subsection{Spectral Radius}
\begin{lemma}
\label{lem:D*}
Let $r$ be such that $\|\Lp^n\|\leq Cr^n$.
Then $D_n\leq \delta^{\alpha-1}Cr^n$.
\end{lemma}
\begin{proof}
For each $W\in \Sigma$,
\[
\delta^{\alpha-1}|W|^{-\alpha}\int_W|DT^{-k}|dm=\delta^{\alpha-1}|W|^{-\alpha}\int_W\Lp^k1
\leq \delta^{\alpha-1}\|\Lp^k1\|_s\leq \delta^{\alpha-1}Cr^k.
\]
Taking the supremum over $W$ and $0 \leq k \leq n$ yields the lemma.
\end{proof}

\begin{rem}
Lemma~\ref{lem:D*} implies that $D_*$ is bounded by the spectral radius of $\Lp$ and since the Lasota-Yorke
estimates imply the reverse inequality, we conclude that in fact $D_*$ is the spectral radius of $\Lp$ on $\B$.
\end{rem}

Thanks to Lemma~\ref{lem:distribution} and Lemma~\ref{lem:D*} we can prove the following
characterization.

\begin{lemma}
\label{lem:spectra}
The spectral radius of $\Lp$ on $\B$ is one and the essential spectral
radius is $\tau$. In addition, calling $\V$
the eigenspace associated to the eigenvalues of modulus one, then
$\Lp$ restricted to $\V$ has a semi-simple spectrum (no Jordan
blocks). Finally, $\V$ consists of signed measures. 
\end{lemma}
\begin{proof}
Recall that by quasi-compactness, the part of the
spectrum larger than $\tau D_*$ is of finite rank (see \cite{Bal}). Now, let
$z$ be in the spectrum of $\Lp$, $|z| > \max\{1,\tau D_*\}$.  Then there must exist an
$h\in\B$ such that $\Lp h=z h$. Accordingly, for each $\vf\in\C^1$, since 
$\vf\circ T^n\in\C^1_{S_n^+}$ for all $n\in\N$,   
\[
|h(\vf)|=|z|^{-n}|\Lp^n h(\vf)|\leq |z|^{-n}|h(\vf\circ T^n)|\leq
|z|^{-n}C \|h\|(|\vf|_\infty+\mu_+^n|D^s\vf|_\infty)
\]
by Lemma~\ref{lem:distribution}.

Thus, if $|z|>1$, we have $h(\vf)=0$ for each $\vf\in\C^1$, which
implies $h=0$ by Remark \ref{rem:inclusion}. 
On the other hand, if $|z|=1$, then it follows that 
$|h(\vf)|\leq C \|h\|\cdot|\vf|_\infty$, so $h$ is a measure.

Next, suppose $\tau D_*\geq 1$.  By the preceding paragraph, the spectral radius of $\Lp$ can be
at most  $\tau D_*$, thus $D_*\leq \tau D_*$, which is impossible
since $\tau<1$. Hence, the spectral radius is $D_* = 1$.

It remains to show that there are no Jordan blocks corresponding
to the peripheral spectrum.
Indeed, suppose that there exists $z\in\mathbb{C}$ and
$h_0,h_1\in\B$ such that $|z|=1$ and $h_0\neq 0$, $\Lp h_0=zh_0$, $\Lp h_1=z
h_1+h_0$.
This would imply $z^{-n}\Lp^n h_1=nz^{-1}h_0+h_1$, and thus
\[
n|h_0(\vf)|\leq |h_1(\vf)| + C \|h_1\| (|\vf|_\infty+\mu_+^n|D^s\vf|_\infty) .
\]
Dividing by $n$ and taking the limit as $n$ approaches infinity, 
it follows that $h_0=0$, contrary to the hypothesis.
\end{proof}

\begin{rem}
\label{rem:D*} Note that Lemma~\ref{lem:spectra} implies $\|\Lp^n\|\leq
C$ for each $n\in\N$, hence Lemma~\ref{lem:D*} implies $D_n\leq C_\delta$ for
all $n\in\N$.
\end{rem}


\subsection{Peripheral Spectrum}
\label{peripheral}

The following two lemmas prove Theorem~\ref{thm:invariant measures},
points (1-4) and part of (5).  
The rest is proved in Section~\ref{statistics}.

Let $\V_\theta$ be the eigenspace associated to the eigenvalue $e^{2\pi i\theta}$. For the
rest of this section, we use $m$ to denote normalized Riemannian volume 
(Lebesgue measure) on $\M$.

\begin{lemma}
\label{lem:inv-prob} 
Recall that $\bar{\mu}:= \Pi_0 1$.  Then,
\begin{enumerate}
  \item[(i)] All the measures in $\V$ are absolutely continuous 
with respect to $\bar \mu$.   Moreover, $1$ belongs to the spectrum.
  \item[(ii)] There exists a finite number of $q_i\in\N$ such 
that the spectrum on the unit disk is 
$\bigcup_k\{e^{2\pi i\frac p {q_k}}\;:\; 0< p\leq q_k,\;p\in\N\}$.
 In addition, the set of ergodic probability measures absolutely continuous with respect to $\bar \mu$ form a basis of $\V_0$. 
  \item[(iii)] For each
$\mu\in\V$, $n\in\N$, we have $\mu(\Si^\pm_{n,\epsilon}) \leq C_n \epsilon^\alpha$.  In particular, 
$\mu(\Si^\pm_n)=0$. 
\end{enumerate}
\end{lemma}
\begin{proof}
(i)
Let $\Pi_\theta$ be the eigenprojector on $\V_\theta$. 
The fact that the spectrum outside the circle of radius $\tau$ consists of only
finitely many eigenvalues of finite multiplicity implies that the limit
\begin{equation}
\lim_{n\to\infty}\frac 1n\sum_{k=0}^{n-1}e^{-2\pi i\theta k}\Lp^k=\Pi_\theta
\end{equation}
is well-defined in the uniform topology of $L(\B,\B)$. Moreover, $\Pi_0$ is
obviously a positive operator and, by density, $\V_\theta=\Pi_\theta\C^1$. 

Accordingly, for each $\mu\in\V_\theta$, there exists $h\in \C^1$ such that
$\Pi_\theta h=\mu$. Thus, for each $\vf\in\C^1$
\begin{equation}
\label{eq:abs-cont}
|\mu(\vf)|=|\Pi_\theta h(\vf)|\leq |h|_\infty \Pi_0 1(|\vf|)=: |h|_\infty\bar\mu(|\vf|).
\end{equation}
That is, each probability measure $\mu\in\V_\theta$ is absolutely continuous with respect to
$\bar \mu$. Moreover, setting $h_\mu:=\frac{d\mu}{d\bar\mu}$, we have $h_\mu\in L^\infty(\M,\bar\mu)$.
This implies $\bar\mu\neq 0$, otherwise the spectral radius of $\Lp$
would be strictly smaller than one, which, recalling Remark
\ref{rem:inclusion}, yields the contradiction 
\[
1=|m(1)|=|\Lp^nm(1)|=\lim_{n\to\infty} |\Lp^nm(1)|\leq \lim_{n\to\infty} C\|\Lp^nm\|=0.
\]
Hence one belongs to the spectrum.
\medskip
\noindent

(ii) $\qquad$
Next, for $\mu\in\V_\theta$ and each $\vf\in\C^1$
\[
\begin{split}
\int \vf h_\mu d\bar \mu&=\mu(\vf)=e^{-2\pi
i\theta}\Lp\mu(\vf)=e^{-2\pi i\theta}\mu(\vf\circ T) \\
&=e^{-2\pi
i\theta}\int \vf\circ T h_\mu d\bar\mu= e^{-2\pi i\theta}\int \vf
h_\mu\circ T^{-1}d\bar\mu. 
\end{split}
\]
Accordingly $ h_\mu\circ T^{-1}=e^{2\pi i\theta}h_\mu$, $\bar\mu$ a.e. 
This in turn means that, setting,
$h_{\mu,k}:=(h_\mu)^k\in L^\infty(\M,\bar\mu)$, since the measure
$d\mu_k:=h_{\mu,k}d\bar\mu$ belongs to $\B$ for each
$k\in\N$,\footnote{Just consider $h\in\C^1$ such that
$\bar\mu(|h-h_{\mu,k}|)\leq\ve$. Then 
setting $d\nu:=hd\bar\mu$,
\[
(\Pi_{k\theta}\nu-\mu_{k})(\vf)\leq \lim_{n\to\infty}\frac
1n\sum_{j=0}^{n-1}\bar\mu(|h-h_{\mu,k}|\circ T^{-j})|\vf_\infty| 
\leq\ve|\vf|_\infty.
\]
Hence, $\mu_k$ is an accumulation point of elements of $\V$ and so it 
belongs to $\V$.
}
 then $\Lp\mu_k=e^{2\pi ik\theta}\mu_k$. That is, $e^{2\pi ik\theta}$
belongs to the peripheral spectrum and since the peripheral spectrum consists
of a finite number of points, it must be that $\theta\in\mathbb{Q}$.  

Now let $\mu\in\V_0$ and choose $h\in\C^1$ such that $\mu=\Pi_0h$. We can
then write $h=h_+-h_-$, $h_\pm:=\max\{0,\pm h\}$. Since $h_\pm$ are
Lipschitz functions, they belong to $\B$. We 
can then define $\mu_\pm:=\Pi_1 h_\pm$.
Thus $\V_0$ is the span of a convex set of probability measures.  

Next, assume that $\Omega$ is an invariant set of positive $\bar\mu$ measure, then for each $\ve>0$ there exists a smooth 
function $\vf_\ve>0$ such that $\bar\mu(|\Id_\Omega-\vf_\ve|)\leq \ve$. Thus, for each continuous function $\phi$
\[
\bar\mu(\vf_\ve\phi\circ T^n)=\bar\mu(\Id_\Omega\phi\circ T^n)+\Or(\ve|\phi|_\infty)
=\bar\mu(\Id_\Omega\phi)+\Or(\ve|\phi|_\infty)=\bar\mu(\vf_\ve\phi)+\Or(\ve|\phi|_\infty).
\]
On the other hand
\[
\lim_{n\to\infty}\frac 1n\sum_{k=0}^{n-1}\bar\mu(\vf_\ve\phi\circ T^n)=\Pi_1(\bar\mu_\ve)(\phi)
\]
where $\bar\mu_\ve(\phi):=\bar\mu(\vf_\ve\phi)$. By the arbitrariness of $\ve$ it follows, setting 
$\bar\mu_\Omega(\phi):=\bar\mu(\Id_\Omega\phi)$, that $\bar\mu_\Omega\in\V_0$ and, since $\V_0$ is finite dimensional, 
that there are only finitely many invariant sets of positive $\bar\mu$ measure and the ergodic decomposition yields a basis of $\V_0$.

\medskip
\noindent
(iii) $\qquad$
Finally, let $\mu\in\V$. By hypothesis, the tangent space of
$\Si^-_m$ is bounded away from $C^s$.  
Calling $\Si^-_{m,\epsilon}$ an $\epsilon$
neighborhood of $\Si^-_m$, set
$\mu_\epsilon(\vf):=\mu(\Id_{\Si^-_{m,\epsilon}}\vf)$. Let $h_n$ be a
sequence that converges to $\mu$ in $\B$, then it is immediate to
check that $h_{n,\epsilon}(\vf):=h_n(\Id_{\Si^-_{m,\epsilon}}\vf)$
belongs to $\B_w$. In addition, 
\[
\int_W\vf h_{n,\epsilon} dm =\int_{W\cap\Si^-_{m,\epsilon}}\vf h_n dm \leq C_m \|h_n\|\epsilon^{\alpha},
\]
for $\vf \in \C^1(W)$.
In the same way one has that $h_{n,\epsilon}$ is a Cauchy sequence in
$\B_w$, thus it must converge to
$\mu_\epsilon(\vf):=\mu(\Id_{\Si^-_{m,\epsilon}}\vf)$.  Since
$\mu_\epsilon(1)\leq C_m \epsilon^\alpha$, the regularity of
$\mu$ implies $\mu(\Si^-_m)=0$. The result follows since
$T\Si^+=\Si^-$. 
\end{proof}

\begin{rem}\label{rem:manifolds}
For uniformly hyperbolic maps with singularities, the estimate
$\bar{\mu}(\Si^\pm_\epsilon) \leq C\epsilon^\alpha$
is enough to conclude the existence of stable and unstable manifolds
for $\bar{\mu}$-a.e.\ $x \in \M$ using the standard Borel-Cantelli
argument (see, for example \cite{liverani wojtkowski}).
In fact, from the ergodic decomposition proved in Lemma~\ref{lem:physical},
it follows that stable and unstable manifolds exist for $m$-a.e.\
$x \in \M$.
\end{rem}

Recall that to each physical measure $\mu$ is associated a 
positive Lebesgue measure invariant set $B_\mu$ such that, for every continuous
function $f$,
\[
\lim_{n\to\infty}\frac 1n\sum_{i=0}^{n-1}f(T^ix)=\mu(f)\quad \forall x
\in B_\mu.
\]

\begin{lemma}
\label{lem:physical}
Let $T$ be a piecewise uniformly hyperbolic map as described in 
Section~\ref{setting}.
\begin{enumerate}
  \item[(i)] $T$ admits only finitely many physical probability measures and 
they belong to $\V_0$.
  \item[(ii)] The ergodic decomposition with respect to Lebesgue and with respect to $\bar\mu$ coincide.
  \item[(iii)] The forward average for each continuous function is well 
defined $m$-almost everywhere and the ergodic decomposition with respect to Lebesgue corresponds to the supports of the physical measures.
\end{enumerate}
\end{lemma}

\begin{proof}
(i)
Let $\mu$ be a physical measure and take a density point $x$ of the associated
set $B_\mu$. Then for each $\ve>0$ there exists an open set $U$ containing 
$x$ such that $m(B_\mu\cap U)\geq (1-\ve)m(U)$. Consider a smooth
probability measure $\mu_U$ supported in $U$, such that
$\mu_U(B_\mu)\geq 1-2\ve$. Then for each $f\in\C^0$,
\[
\begin{split}
\Pi_1\mu_U(f)=&\lim_{n\to\infty}\frac
1n\sum_{k=0}^{n-1}\mu_U(f\circ T^i)=\lim_{n\to\infty}\frac
1n\sum_{k=0}^{n-1}\mu_U(f\circ T^i\Id_{B_\mu})+{\mathcal
O}(|f|_\infty\ve)\\
& =\mu_U(\Id_{B_\mu})\mu(f)+{\mathcal O}(|f|_\infty\ve)
=\mu(f)+{\mathcal O}(|f|_\infty\ve).
\end{split}
\]
This means that $\mu$ can be approximated by elements of $\V_0$ and
therefore $\mu\in\V_0$. In addition, since the physical 
measures are ergodic by definition, it follows that they must belong 
to the ergodic elements of $\V_0$.

\medskip
\noindent
(ii) $\qquad$
Consider an invariant set $A$ of positive $\bar \mu$ measure such that 
$\bar{\mu}$ restricted to $A$ is ergodic. If we consider the 
set $A_u:=\bigcup_{x\in A}W^u(x)$,\footnote{By $W^u(x),\, W^s(x)$ we designate the unstable and stable manifolds of $x$, which exists $\bar\mu$-a.e. 
(see Remark \ref{rem:manifolds}). Hence, by eventually changing $A$ on 
a set of zero $\bar{\mu}$-measure, we can assume $A\subset A_u$.} 
then $T^{-1}A_u\subset A_u$ and $\bar\mu(A_u\setminus A)=0$. 
Indeed, for each continuous function $\vf$ the backward average on 
$A_u$ coincides with the backward average on $A$ and hence equals 
$\bar\mu(\vf\Id_A)\bar\mu(A)^{-1}$ $\bar\mu$-a.e., by ergodicity. 
In addition, for each $\ve>0$ there exists $\vf\in\C^0$ such that 
$\bar\mu(|\vf-\Id_A|)\leq \ve$.  But then
\[
\lim_{n\to\infty}\frac 1n\sum_{k=0}^{n-1}\vf\circ T^{-k}\Id_{A_u}=\bar\mu(\vf\Id_A)\bar\mu(A)^{-1} \Id_{A_u}\quad \bar\mu\text{-}a.e.
\]
Integrating, we have
$\bar\mu(\vf\Id_{A_u})=\bar\mu(\vf\Id_A)\bar\mu(A)^{-1} \bar\mu(A_u)$, which yields 
\[
\bar\mu(A)=\bar\mu(A_u)+\Or(\ve)
\]
which, by the arbitrariness of $\ve$, proves the claim. 

Next, consider the neighborhood of the singularity $\Si_\delta$. By Lemma \ref{lem:inv-prob}(iii) it follows that 
$\bar\mu(\Si_\delta)\leq C\delta^\alpha$, hence setting 
$\Si_{\delta}^*:=\bigcup_{m\in\mathbb{Z}}T^{m}\Si_{\delta |m|^{-2/\alpha}}$ we have 
$\bar\mu(\Si_\delta^*)\leq C\delta^\alpha\sum_{m\in\mathbb{Z}} |m|^{-2}\leq C\delta^\alpha$.
Hence, by choosing $\delta$ small enough and setting $B:=A_u\setminus \Si_\delta^*$, it holds that $\bar\mu(B)>0$. 
In addition, for each $x\in B$ the manifolds $W^u(x), W^s(x)$ have length at least $\delta$ (see, e.g., \cite{liverani wojtkowski}).  
Define $B_\epsilon$ to be the $\epsilon$-neighborhood of $B$.

Finally, for each $n\in\N$ let 
$U_n=\bigcup_{\substack{K\in{\mathcal K}_n\\ K\cap B\neq \emptyset}} K\bigcap B_{\delta/2}$, and $\tilde U_n:= U_n\bigcap B_{\delta/4}$,
clearly $U_{n+1}\subseteq U_n$. 
Setting $\tilde B:=\bigcap_{n\in\N}U_n$, we have 
$\bar\mu(\tilde B \setminus A_u)=0$. Indeed, if $x\in\tilde B\setminus A_u$, 
then, for each $n\in\N$, there exist $y_n\in A_u\setminus \Si_\delta^*$ and $K_n\in\K_n$ such that 
$x,y_n\in K_n$. On the other hand each $W^s(y_n)$ has size $\delta$, by construction. We can then extract a 
converging subsequence $\{y_{n_j}\}$, and call $W$ the limit of $W^s(y_n)$. It is easy to check that $W$ is a stable 
manifold and since in the transversal direction the map is expanding it must be $W=\bigcap_n K_n$. 
Accordingly, $W\subset W^s(x)$. The uniform transversality 
of the stable and unstable manifolds implies then that 
$W^u(y_n)\cap W^s(x)\neq \emptyset$ for $n$ large enough. 
Hence, $W^s(x)\cap A_u\neq \emptyset$, that is $\tilde{B} \subset A_{u,s} := \bigcup_{x \in A_u} W^s(x)$. 
The claim follows then since one can prove $\bar\mu(A_{u,s}\setminus A_u)=0$ by 
the same argument that established $\bar\mu(A_u\setminus A)=0$.

In addition, by the regularity of both $\bar{\mu}$ and $m$, 
for each $k\in\N$ and $\ve>0$, there exists $n_\ve$ such that
\[
\bar\mu(U_{n_\ve}\setminus \tilde{B})+\sum_{j=0}^k\Lp^km(U_{n_\ve} \setminus \tilde{B})\leq \ve.
\]
Let $\phi$ be a smooth function such that $\phi_{\upharpoonright_{B_{\delta/4}}}=1$ and 
$\text{supp}\,\phi$ is contianed in $B_{\delta/2}$.  Clearly $\phi \Id_{U_n}\in\C^1_{\Si_n^+}$,
we can the use the above inequality and Lemma \ref{lem:distribution} to write
\[
\begin{split}
&\bar\mu(\tilde B)=\bar\mu(\tilde U_{n_\ve})+\Or(\ve)\leq \bar\mu(\Id_{U_{n_\ve}}\phi)+\Or(\ve)
=\frac 1k\sum_{j=0}^{l-1}\Lp^jm(\Id_{U_{n_{\ve}}}\phi)+\Or(k^{-1}+\ve)\\
&=\frac 1k\sum_{j=0}^{l-1}\Lp^jm(\tilde B)+\Or(k^{-1}+\ve)
\leq \frac 1k\sum_{j=0}^{l-1}\Lp^jm(A_{u,s})+\Or(k^{-1}+\ve)=m(A_{u,s})+\Or(k^{-1}+\ve).
\end{split}
\]
By choosing first $k$ large enough and then $\ve$ sufficiently small, we have $m(A_{u,s})>0$.

Next notice that if $A_{u,s}$ is not an ergodic component for $m$, then there exists an invariant 
set of positive $m$-measure that will support a physical measure, but such a measure would belong 
to $\V_0$ by point (i) and this would be a contradiction.

\medskip
\noindent
(iii) $\qquad$
The preceding argument also implies that if $\mu_i$ is a basis of $\V_0$ made of ergodic measures, 
then they are physical measures and $\{B_{\mu_i}\}$ corresponds 
to the ergodic decomposition with respect to Lebesgue. In addition, since 
$m(B_{\mu_i}) \geq \bar{\mu}(B_{\mu_i})$ and $\sum_i \bar{\mu}(B_{\mu_i}) =1$, the inequality must be
an equality and the forward average for each continuous function is well 
defined $m$-almost everywhere.
\end{proof}


\subsection{Statistical Properties and Ruelle Resonances}
\label{statistics}

In addition to providing information about the invariant measures, the established
spectral picture has other far reaching implications. To discuss them
let us define the {\em correlation functions}. For each
$f,g\in\C^\beta$ define
\[
C_{f,g}(n):=\bar\mu(f g\circ T^n)-\bar\mu(f)\bar\mu(g).
\]
If the system is mixing (that is, one is the only eigenvalue on
the unit circle and it is simple), then for each $\sigma$ larger than the norm of the
second largest eigenvalue (or $\tau$ if no other eigenvalue is present
outside the essential spectral radius) holds
\begin{equation}
\label{eq:decay}
|C_{f,g}(n)|\leq C\sigma^n|f|_{\C^\beta}|g|_{\C^\beta}.
\end{equation}
In other words we have the well-known dichotomy: either the system does
not mix or it mixes exponentially fast (on H\"older observables).

More generally, we can define the Fourier transform of the correlation
function:
\[
\hat C_{f,g}(z):=\sum_{n\in\mathbb Z}z^nC_{f,g}(n).
\]
The above quantity is widely used in the physics literature where usually
one assumes that it is convergent in a neighborhood of $|z|=1$ (here
this follows already from \eqref{eq:decay}) and it
has a meromorphic extension on some larger annulus. The poles of such a
quantity are, in principle, measurable in a physical system and are
called {\em Ruelle resonances} 
(see \cite{ruelle 1, ruelle 2, parry poll 1, parry poll 2}). 
Due to our results we can substantiate
the above picture for the class of systems at hand.

Indeed, note that we can assume, without loss of generality,
$\bar\mu(f)=\bar\mu(g)=0$ and that if we define $\mu_f(\vf):=\bar\mu(f\vf)$,
$\mu_g(\vf):=\bar\mu(g\vf)$, then $\mu_f,\mu_g\in\B$. thus,
\[
\begin{split}
\hat C_{f,g}(z)&=\sum_{n=0}^\infty z^n\bar\mu(f g\circ
T^n)+\sum_{n=0}^\infty z^{-n} \bar\mu(f\circ T^{n} g) -\bar\mu(fg)\\
&=\sum_{n=0}^\infty z^n\Lp^n\mu_f(g)+\sum_{n=0}^\infty z^{-n}
\Lp^n\mu_g(f) -\bar\mu(fg)\\
&=(z-\Lp)^{-1}\mu_f(g)+(z^{-1}-\Lp)^{-1}\mu_g(f)-\bar\mu(fg).
\end{split}
\]
It is thus obvious that the desired meromorphic extension is provided
by the resolvent and that the poles are in one-to-one correspondence
(including multiplicity) with the spectrum of $\Lp$. More precisely we
have a meromorphic extension in the annulus $\{z\in{\mathbb C}\;:\;
\tau<|z|<\tau^{-1}\}$.
\begin{rem}
Note that the above fact shows that the spectral data of the operator
$\Lp$ on $\B$ is not a mathematical artifact but has a well-defined
meaning which does not depend on any of the many arbitrary choices we
have made in the construction of our functional analytic setting.
\end{rem} 
\begin{rem}
\label{rem:optimal}
In the present situation the best one can do is to choose
$\alpha=\beta=q=\frac 12$; moreover, if one assumes that $M(n)$ grows
sub-exponentially (this is the case for billiards), then one has
(assuming for simplicity $\lambda^{-1}=\mu_+$) that $\tau$ can be
chosen arbitrarily close to $\lambda^{-\frac 12}$. At the moment it is unclear if
such an estimate for the size of the meromorphic extension is real or is
an artifact of the method of proof.
\end{rem}
Another result that can be easily obtained by the present method is
the Central Limit Theorem. Let $f\in\C^\beta$ with
$\bar\mu(f)=0$ and define $S_n(f):=\sum_{k=0}^{n-1}f\circ T^k$. Then
\[
\bar\mu(e^{-izS_n})=\Lp_z^n\bar\mu(1)
\]
where $\Lp_z$ is the operator defined by
$\Lp_zh(\vf):=h(e^{-izf}\vf\circ T)$. Since $\Lp_z$ depends
analytically on $z$, one can use standard perturbation theory to show
that the leading eigenvalue is given by $1-\sigma z^2$, where $\sigma$
is the variance. Accordingly
\[
\lim_{n\to\infty}\bar\mu(e^{-i\frac{z}{\sqrt n}S_n})=\lim_{n\to\infty}
\left(1-\frac{\sigma z^2}n\right)^n=e^{-\sigma z^2}
\]
which is exactly the CLT. Other types of results (e.g.\ large
deviations) can be approached along similar lines (see
\cite{hennion, chazottes} for more details).


\section{Perturbation Results}
\label{perturbations}

Recall from Section~\ref{pert results} the set $\Gamma_{B_*}$ of maps $\tT$ that satisfy
the same assumptions as $T$ in Section~\ref{setting}. 
In this section we derive results for several classes of perturbations
and prove Theorems~\ref{thm:pert} and \ref{theorem:holes}.


\subsection{Deterministic Perturbations}
\label{sec:det-pert}

\begin{lemma}
\label{lem:perturbation}
If two maps $T_1,\,T_2 \in \Gamma_{B_*}$ satisfy $\gamma (T_1,T_2) \leq \ve \leq \ve_0$, 
then for each $h \in\B$,
\[
|\Lp_{T_1} h-\Lp_{T_2} h|_w \leq C_b \ve^\beta \|h\|.
\]
\end{lemma}
\begin{proof}
For $\ve \leq \ve_0$, we may choose the set of approximate stable leaves
$\Sigma$ so that $T^{-1}_i \Sigma \subset \Sigma$ for $i=1,2$.  And similarly
for the approximate unstable family $\F^u$.

We first fix a leaf $W \in \Sigma$ and $\vf$ with $|\vf|_{\C^1(W)} \leq 1$ and write
\[
\int_W (\Lp_{T_1} - \Lp_{T_2})h \, \vf \, dm =
\int_{T_1^{-1}W}h |DT_1|^{-1} J_WT_1 \vf \circ T_1 - \int_{T_2^{-1}W}h|DT_2|^{-1} J_WT_2 \vf \circ T_2.
\]
Away from singularities, $T_1^{-1}W$ and $T_2^{-1}W$ are $\ve$-close so we may partition $T_1^{-1}W$
and $T_2^{-1}W$ as we did in Section~\ref{unstable norm}.

Let $N_\ve^-$ denote the $\ve$ neighborhood of the union of the singularity curves of $T_1^{-1}$ 
and $T_2^{-1}$.  Consider one component $U_j$ of $W\backslash N_\ve^-$.  By assumption, 
we may choose functions
$F^i_j$ defining the curves $T_i^{-1}U_j$ such that $d_\Sigma (T_1^{-1}U_j, T_2^{-1}U_j) \leq \ve$.
(If $\max \{ |T_1^{-1}U_j|, |T_2^{-1}U_j| \} > 2 \delta$, we further subdivide
$U_j$ so that all components of $T_1^{-1}U_j$ and $T_2^{-1}U_j$ have length between
$\delta$ and $2\delta$.)

Denote by $V_j$ the connected components of $W \cap N_\ve^-$ and note that $|V_j| \leq C\ve$
and that there are at most $L+2$ such pieces.

We estimate the integrals over the pieces $T_i^{-1}V_j$ similarly to \eqref{eq:first unstable}
\begin{equation}
\label{eq:short pieces close}
\sum_{i,j} \int_{T_i^{-1}V_j} h |DT_i|^{-1} J_WT_i \, \vf \circ T_i \, dm 
\leq C \|h\|_s \sum_{i,j} |V_j|^\alpha \lambda^{-1} \mu^{-\alpha}
\leq C \|h\|_s \ve^\alpha   .
\end{equation}

We split up the integrals over the $T_i^{-1}U_j$ as follows,
\begin{eqnarray}
\sum_j & \! \! \! \! \! & \! \! \! \! \!
\int_{T_1^{-1}U_j} h |DT_1|^{-1} J_WT_1 \, \vf \circ T_1 \, dm -
\int_{T_2^{-1}U_j} h |DT_2|^{-1} J_WT_2 \, \vf\circ T_2  \, dm  \nn \\ 
 & = & \sum_j \int_{T_1^{-1}U_j} h |DT_1|^{-1} J_WT_1 \, \vf \circ T_1 \, dm - 
                  \int_{T_2^{-1}U_j} h f \, dm  \label{eq:long pieces close} \\
 &   & + \sum_j \int_{T_2^{-1}U_j} h (f- |DT_2|^{-1}J_WT_2 \vf\circ T_2) dm  \nn
\end{eqnarray}
where $f = [ |DT_1|^{-1} J_WT_1 \, \vf\circ T_1 ] \circ G_{F_j^1} \circ G_{F_j^2}^{-1}$.
Note that
$d_q(|DT_1|^{-1} J_WT_1 \, \vf \circ T_1, f) = 0$ so that the first term of 
\eqref{eq:long pieces close} can be estimated by
\begin{equation}
\label{eq:u-pert}
\sum_j \left| \int_{T_1^{-1}U_j} h |DT_1|^{-1} J_WT_1 \vf \circ T_1 - 
                  \int_{T_2^{-1}U_j} h f \right|
\leq C\ve^\beta \|h\|_u .
\end{equation}

We estimate the second term of \eqref{eq:long pieces close} using the strong stable norm.
We follow \eqref{eq:diff} to estimate the $\C^q$-norm of the functions involved.
\[
\begin{split}
& |f- |DT_2|^{-1}J_WT_2 \vf\circ T_2 |_{\C^q(T_2^{-1}U_j)} \\
& \leq C |[ |DT_1|^{-1} J_WT_1 \vf\circ T_1 ] \circ G_{F_j^1}  
 - [ |DT_2|^{-1}J_WT_2 \vf\circ T_2 ] \circ G_{F_2^2}|_{\C^q(I_{r_j})} \\
& \leq C |\vf\circ T_1 \circ G_{F_j^1} - \vf\circ T_2 \circ G_{F_2^2}|_{\C^q(I_{r_j})} \\
& \qquad + C | (|DT_1|^{-1} J_WT_1 ) \circ G_{F_j^1} 
   - (|DT_2|^{-1}J_WT_2) \circ G_{F_2^2}|_{\C^q(I_{r_j})}.
\end{split}
\]
The first term can be bounded using an estimate analogous to \eqref{eq:vf-unstable}
and Lemma~\ref{lem:change leaf}.  The second term can be bounded using Lemma~\ref{lem:distortion} 
and the fact that $|T_1 -T_2|_{\C^2} < \ve$ on
$U^i_j$.  Putting these estimates together, we
conclude that $|f- |DT_2|^{-1}J_WT_2 \vf\circ T_2 |_{C^q} \leq C \ve^{1-q}$ so we may
estimate the second term of \eqref{eq:long pieces close} by
\[
\int_{T_2^{-1}U_j} h (f - |DT_2|^{-1}J_WT_2 \vf\circ T_2) \; \leq \; C \ve^{1-q} \|h\|_s.
\]

Putting this estimate together with \eqref{eq:short pieces close} and \eqref{eq:u-pert}, we have
\begin{equation}
\label{eq:operator close}
\left| \int_W \Lp_{T_1} h \vf dm - \int_W \Lp_{T_2} h \vf dm \right|
\leq C (\|h\|_s \ve^\alpha + \|h\|_u \ve^\beta + \|h\|_s \ve^{1-q}) \leq Cb^{-1} \ve^\beta \|h\|.
\end{equation}
Taking the supremum over all $W \in \Sigma$ and $\vf \in \C^1(W)$ yields the lemma.
\end{proof}

Lemma~\ref{lem:perturbation} implies $||| \Lp_{T_1} - \Lp_{T_2} ||| \leq C \ve^\beta$
whenever $\gamma(T_1, T_2) \leq \ve$.  Since both $T_1$ and $T_2$ satisfy the Lasota-Yorke
inequalities \eqref{eq:weak norm}-\eqref{eq:unstable norm}, we may
apply the results of \cite{KL} to our operator $\Lp:\B \to \B_w$.


\subsection{Smooth Random Perturbations}
\label{sec:random-pert}

Recall the transfer operator $\Lp_{\nu,g}$ associated with the random process defined in
Section~\ref{results}.  For the remainder of this section, we fix constants $\lambda$, $\mu$,
$\mu_+$ and $D_n$
such that \eqref{eq:exp def} and \eqref{eq:def-Dn} are satisfied for all $\tT \in X_\ve$.

The following is a generalization of Lemma~\ref{lem:perturbation} which shows that the transfer
operator associated with the random perturbation is also close to $\Lp_T$ in the sense of \cite{KL}.

\begin{lemma}
\label{lem:random perturbation}
$ \displaystyle ||| \Lp_{\nu,g} - \Lp_T ||| \leq C_b A \ve^\beta$.
\end{lemma}
\begin{proof}
Let $h$, $\vf \in \C^1(\M)$, $|\vf|_{\C^1} \leq 1$, and $W \in \Sigma$.  Then
using \eqref{eq:operator close} of Lemma~\ref{lem:perturbation},
\begin{eqnarray*}
\left| \int_W \Lp_{\nu,g} h \, \vf \; dm - \int_W \Lp_T h \, \vf \; dm \right| 
   & = & \left| \int_\Omega \int_W (\Lp_{T_\omega} h(x)  - \Lp_T h(x)) 
     \, \vf(x) \, g(\omega, T_\omega^{-1}x) \; dm d\nu  \right| \\
   & \leq & \int_\Omega C_b \ve^\beta \|h\| |g(\omega, \cdot)|_{C^1} d\nu(\omega)
   \; \; \leq \; \; C_b A \ve^\beta \|h\|.
\end{eqnarray*}
\end{proof}

We next prove uniform Lasota-Yorke
estimates for the operator $\Lp_{\nu,g}$.
First, we need to introduce some notation.  Let
$\ob_n = (\omega_1, \ldots ,\omega_n) \in \Omega^n$.  We define 
$T_{\ob_n} = T_{\omega_n} \circ \cdots \circ T_{\omega_1}$ and similarly
$DT_{\ob_n} = \Pi_{j=1}^{n} DT_{\omega_j}(T_{\ob_{j-1}})$.

\begin{lemma}
\label{lem:uniform estimates}
Let $\Delta(\nu,g) \leq \ve$.
For $\ve$ sufficiently small,
there exists $\delta_0>0$ and a constant $C = C_{a,A}$, such that for all $h \in \B$, 
$\delta \leq \delta_0$ and $n \geq 0$, 
$\Lp_{\nu,g}$ satisfies
\begin{eqnarray*}
| \Lp^n_{\nu,g} h |_w   & \leq & C D^n |h|_w \; ,  \\
\| \Lp_{\nu,g}^n h\|_s & \leq & C \max\{ \rho, \mu_+^q\}^n D_n \|h\|_s + C_\delta D_n |h|_w \; ,  \\
\| \Lp^n_{\nu,g} h\|_u & \leq & C \lambda^{-\beta n} D_n \|h\|_u 
       + C (D_n + L_n \lambda^{-n} \mu^{-\alpha n}) \|h\|_s \; .
\end{eqnarray*}
\end{lemma}
\begin{proof}
The proofs follow from those of Section~\ref{inequalities}, 
except that we have the added function $g(\omega,x)$.
Notice that 
\[
\Lp_{\nu,g}^n h(x) = \int_{\Omega^n} h \circ T_{\ob_n}^{-1} |DT_{\ob_n}(T_{\ob_n}^{-1})|^{-1} 
          \Pi_{j=1}^n g(\omega_j, T_{\omega_j}^{-1} \circ \cdots 
          \circ T_{\omega_n}^{-1}x) \; d\nu^n(\ob_n)  .
\]
{\bf Estimating the strong stable norm.}  For any $W \in \Sigma$, 
we define the connected pieces $W_i$ of $T_{\ob_n}^{-1}W$ inductively just as we did for $T^{-n}W$ in 
Section~\ref{weak norm}.  Following the estimates of Section~\ref{stable norm}, we write
\begin{equation}
\label{eq:stable random split}
\begin{split}
\int_W \Lp_{\nu,g}^n h & \, \vf \; dm =
               \int_{\Omega^n} \sum_i \left\{ \int_{W_i} h \overline{\vf}_i |DT_{\ob_n}|^{-1}
               J_WT_{\ob_n} \Pi_{j=1}^n g(\omega_j, T_{\ob_{j-1}} x) \; dm(x)   \right.         \\
        & \left. + \frac{1}{|W_i|} \int_{W_i} \vf \circ T_{\ob_n} \int_{W_i} h |DT_{\ob_n}|^{-1}
               J_WT_{\ob_n} \Pi_{j=1}^n g(\omega_j, T_{\ob_{j-1}} x) \; dm(x) \right\} d\nu^n(\ob_n)
\end{split}
\end{equation}
where $\overline{\vf}_i =  \vf \circ T_{\ob_n} - \frac{1}{|W_i|} \int_{W_i} \vf \circ T_{\ob_n}$.
We fix $\ob_n$ and define 
$G_{\ob_n}(x) = \Pi_{j=1}^n g(\omega_j, T_{\ob_{j-1}} x)$.

To estimate the first sum in \eqref{eq:stable random split}, we note that \eqref{eq:C^q small} implies
\[
|\overline{\vf}_i|_{C^q(W_i)} \leq C |J_WT_{\ob_n}|_{C^0(W_i)}^q |W|^{-\alpha} .
\]
Then, using \eqref{eq:first stable}, we estimate
\begin{equation}
\label{eq:first random}
\begin{split}
\sum_{i} & \int_{W_i} h \; \overline{\vf}_i |DT_{\ob_n}|^{-1} J_WT_{\ob_n} G_{\ob_n} \; dm  \\
   & \leq \sum_i C \|h\|_s |W_i|^\alpha |\, |DT_{\ob_n}|^{-1} J_WT_{\ob_n}|_{C^q(W_i)} 
                                       |\overline{\vf}_i|_{C^q(W_i)} |G_{\ob_n}|_{C^q(W_i)} \\
   & \leq \sum_{i} C \|h\|_s |W_i|^\alpha  
   ||DT_{\ob_n}|^{-1} J_WT_{\ob_n}|_{C^0(W_i)} |J_WT_{\ob_n}|_{C^0(W_i)}^q |W|^{-\alpha}  |G_{\ob_n}|_{\C^q(W_i)}.
\end{split}
\end{equation}
The only additional term here is $|G_{\ob_n}|_{\C^q(W_i)}$, which we now show is bounded independently of $n$ and $W_i$.

\begin{sublem}
\label{lem:G dist}
Let $W_i \in \Sigma$ be a smooth component of $T_{\ob_n}^{-1}W$.  There exists a constant $C>0$, 
independent
of $W$, $n$ and $\ob_n$ such that
\[
\left| \Pi_{j=1}^n g(\omega_j, T_{\ob_{j-1}} \cdot ) \right|_{\C^1(W_i)} 
\leq C \Pi_{j=1}^n g(\omega_j, T_{\ob_{j-1}} x )
\]
for any $x \in W_i$.
\end{sublem}
\begin{proof}
The proof follows the usual distortion estimates along stable leaves.  For any $x$, $y \in W_i$,
\begin{eqnarray*}
\log \frac{\Pi_{j=1}^n g(\omega_j, T_{\ob_{j-1}} x )}{\Pi_{j=1}^n g(\omega_j, T_{\ob_{j-1}} y )}
     & \leq & \sum_{j=1}^n a^{-1} |g(\omega_j, \cdot)|_{C^1(W_i)} d(T_{\ob_{j-1}}x, T_{\ob_{j-1}}y)  \\
     & \leq & \sum_{j=1}^\infty Aa^{-1} C \mu_+^{j-1} d(x,y) \; \; =: \; \; c_0 d(x,y) ,
\end{eqnarray*}
using property (iii) of $g$.
The distortion bound yields the lemma with $C = c_0 e^{c_0}$.
\end{proof}

The sublemma allows us to estimate \eqref{eq:first random} using 
\eqref{eq:first stable}.
\begin{equation}
\label{eq:second random}
\sum_{i} \int_{W_i} h \overline{\vf}_i |DT_{\ob_n}|^{-1} J_WT_{\ob_n} G_{\ob_n} \; dm 
    \leq C \|h\|_s D_n \mu_+^{q n} \Pi_{j=1}^n g(\omega_j, T_{\ob_{j-1}} x_* )
\end{equation}
where $x_*$ is some point in $T_{\ob_n}^{-1}W$.

We estimate the second term of \eqref{eq:stable random split} in a similar way according to
\eqref{eq:second stable}.  Each time, we replace
$|G_{\ob_n}|_{\C^q}$ or $|G_{\ob_n}|_{\C^1}$ according to Sublemma~\ref{lem:G dist}.  
\[
\begin{split}
\sum_i \frac{1}{|W_i|} \int_{W_i} \vf \circ T_{\ob_n} \, dm &  \int_{W_i} h |DT_{\ob_n}|^{-1}
               J_WT_{\ob_n} G_{\ob_n} \, dm \\
& \leq (C\|h\|_s \rho^n + C_\delta D_n |h|_w) \Pi_{j=1}^n g(\omega_j, T_{\ob_{j-1}} x_* )
\end{split}
\] 
Combining this estimate with \eqref{eq:second random}, we have
\[
\int_W \Lp_{T_{\ob_n}}h \vf \; dm \leq (C \|h\|_s (D_n \mu_+^{q n}+ \rho^n) +  C_\delta D_n |h|_w) 
              \Pi_{j=1}^n g(\omega_j, T_{\ob_{j-1}} x_* ) 
\]
Now integrating this expression over $\Omega^n$, we integrate one $\omega_j$ at a time starting 
with $\omega_n$.
Note that $\int_\Omega g(\omega_n, T_{\ob_{n-1}} x_* ) d\nu(\omega_n) = 1$ by assumption on 
$g$ since $T_{\ob_{n-1}} x_*$
is independent of $\omega_n$.  Similarly, each factor in $G_{\ob_n}$ integrates to 1 so that
\[
\| \Lp^n_{\nu, g} h \|_s \leq C \|h\|_s (D_n \mu_+^{q n}+ \rho^n) +  C_\delta D_n |h|_w 
\]
which is the Lasota-Yorke inequality for the strong stable norm.

The inequalities for the strong unstable norm and for the weak norm follow almost identically, always
using Sublemma~\ref{lem:G dist}.
\end{proof}


\subsection{Hyperbolic Systems with Holes}
\label{holes}

We adopt the notation and conditions introduced in Section~\ref{hole results}.
The first lemma shows that we can make the operators $\Lp$ and $\Lp_H$ arbitrarily
close by controlling the ``diameter'' $r$ of the hole along elements of $\Sigma$ 
and the number $P$ of connected components of the hole that a leaf can intersect
at time 1.

\begin{lemma}
\label{lemma:hole close}
Let $H$ be a hole satisfying assumption (H1).
There exists $C>0$ depending only on $T$ such that
\[
||| \Lp - \Lp_H ||| \leq C P r^\alpha .
\]
\end{lemma}

\begin{proof}
Let $h \in \C^1(\M)$, $W \in \Sigma$ and $\vf \in \C^1(W)$ with $|\vf|_{\C^1(W)} \leq 1$.
Recall that $\M^1 \subset \M \backslash H$ is the set of points which remains
in $\M$ until at least time 1.  Let $1_{\M \backslash \M^1}$ denote the indicator
function of $\M \backslash \M^1$. 
\[
\begin{split}
& \int_W (\Lp - \Lp_H)h \; \vf \; dm \; = \; \int_W \Lp (1_{\M \backslash \M^1}h)\; \vf \; dm  \\
   & = \; \int_{T^{-1}W \cap \M \backslash \M^1} h \; \vf\circ T |DT|^{-1} J_WT \; dm
   \; \leq \; \sum_{\tilde{W}_i} \|h\|_s |\tilde{W}_i|^\alpha |\vf \circ T|_{\C^q(\tilde{W}_i)} 
   |\,|DT|^{-1} J_WT|_{\C^0(\tilde{W}_i)} 
\end{split}
\]
where $\tilde{W}_i$ are the connected components of $T^{-1}W \cap \M \backslash \M^1$, i.e. the pieces of
$T^{-1}W$ which are in the hole at time 0 or 1.  
We recall from the estimates of Section~\ref{inequalities}
that $|\vf \circ T|_{\C^q(\tilde{W}_i)} \leq |\vf|_{\C^q(W)}$.  Also, the distortion bound \eqref{eq:distortion}
implies $|J_WT| |\tilde{W}_i| \leq C |T\tilde{W}_i|$.  We then have
\[
\int_W (\Lp - \Lp_H)h \; \vf \; dm \leq C \|h\|_s \sum_i |T\tilde{W}_i|^\alpha \leq C \|h\|_s Pr^\alpha ,
\]
which completes the proof of the lemma. 
\end{proof}

The next proposition proves uniform Lasota-Yorke estimates for $\Lp_H$ which are independent of $H$ satisfying
assumptions (H1) and (H2).
Once it is proven, we may use it in combination with Lemma~\ref{lemma:hole close} to invoke
the results of \cite{KL} and conclude that the spectra of $\Lp$ and $\Lp_H$ are close if $r$ is small.
This proves Theorem~\ref{theorem:holes}.
\begin{proposition}
\label{prop:hole compact}
Let $H$ be a hole satisfying assumptions (H1) and (H2) of Section~\ref{hole results} and let
$\rho_1 := \frac{L+P}{\lambda \mu^\alpha} < 1$.  
Choose $\beta \leq \alpha/2$.
There exists $\delta_0>0$, depending only on $P$, such that for all 
$h \in \B$, $\delta \leq \delta_0$
and $n \geq 0$, $\Lp_H$ satisfies
\begin{eqnarray}
| \Lp^n_H h |_w   & \leq & C D_n |h|_w \; ,   \label{eq:w hole} \\
\| \Lp_H^n h\|_s & \leq & C \max\{ \rho_1, \mu_+^q\}^n D_n \|h\|_s 
               + C_\delta D_n |h|_w \; , \label{eq:ss hole}  \\
\| \Lp^n_H h\|_u & \leq & C \lambda^{-\beta n} D_n \|h\|_u 
         + C( D_n + (L_n + P_n)\lambda^{-n} \mu^{-\alpha n}) \|h\|_s \; .\label{eq:su hole}
\end{eqnarray}
\end{proposition}

\begin{proof}
Our estimates follow closely those of Section~\ref{inequalities}, so
to avoid repetition we indicate only where the presence of the holes requires us to modify those estimates. 
First notice that Lemmas~\ref{lem:counting0} and \ref{lem:counting} hold for the map with holes
with $\rho_1$ in place of $\rho$.  This is because
the definition of the elements $W_i^k$ of $\W_k$ and their tree-like structure remains unchanged.  
The number of connected components of $\tT^{-n}W$ may be greater, but the growth of the number of 
short pieces is controlled
by assumption (H2).  Summing up to most recent long ancestors as we did in the proof of
Lemma~\ref{lem:counting0} and using (H2), we see that equation~(\ref{eq:short pieces}) becomes
\[
\sum_{i \in J_n(W_j^k)} |W_i^n|^\varsigma ||DT^n|^{-1}J_WT^n|_{\C^0(W_i^n)} 
\leq C ||DT^k|^{-1}J_WT^k|_{\C^0(W^k_j)} |W^k_j|^\varsigma \rho_1^{n-k}.
\]
The proof of expressions analogous to equations~\eqref{eq:long remaining}-\eqref{eq:summable} 
is now identical to the proof of Lemma~\ref{lem:counting0}.  We conclude that
\begin{equation}
\label{eq:summable hole}
\sum_i |W^n_i|^\varsigma ||DT^n|^{-1} J_WT^n|_{\C^0(W^n_i)} 
\; \leq \; C D_n \delta^{\varsigma - \alpha} |W|^\alpha 
  + C|W|^\varsigma \rho_1^n.
\end{equation}

\medskip
\noindent
{\bf Estimating the weak norm.}  
For any $h \in \C^1(\M)$, $W \in \Sigma$ and $\vf \in \C^1(W)$ with $|\vf|_{C^1(W)} \leq 1$, we have
\[
\begin{split}
\int_W \Lp_H^n h \, \vf \, dm & \; = \; \sum_{W_i \in \W_n} \int_{W_i} h |DT^n|^{-1} J_WT^n \, \vf\circ T^n \, dm \\
 & \leq \; C |h|_w \sum_{W_i \in \W_n} ||DT^n|^{-1} J_WT^n|_{\C^0(W^n_i)} \; \leq \; C D_n |h|_w
\end{split}  
\]
where in the last inequality we have used \eqref{eq:summable hole} with $\varsigma=0$.  
This proves \eqref{eq:w hole}.

\medskip
\noindent
{\bf Estimating the strong stable norm.}
As in Section~\ref{stable norm}, we define 
$\overline{\vf}_i =  \vf \circ T^n - \frac{1}{|W_i|} \int_{W_i} \vf \circ T^n$.
Equation~\eqref{eq:first stable} remains unchanged,
\[
\sum_i \int_{W_i} h |DT^n|^{-1} J_WT^n \, \overline{\vf}_i \, dm 
\leq C \|h\|_s D_n \mu_+^{q n} .
\]
The estimate for equation~\eqref{eq:second stable} is modified slightly according to \eqref{eq:summable hole},
\[
\sum_i \frac{1}{|W_i|} \int_{W_i} \vf \circ T^n \, dm \int_{W_i} h |DT^n|^{-1} J_WT^n \, dm
\leq C \|h\|_s \rho_1^n + C_\delta D_n |h|_w.
\]
Combining these two estimates, we see that
\[
\| \Lp_H^n h \|_s \leq C \|h\|_s (D_n \mu_+^{q n} + \rho_1^n) + C_\delta D_n |h|_w  ,
\]
which proves~\eqref{eq:ss hole}. 

\medskip
\noindent
{\bf Estimating the strong unstable norm.}
Given two admissible leaves $W^1$ and $W^2$ satisfying $d_\Sigma (W^1, W^2) \leq \ve$, we partition them into
long pieces $U^i_j$ and short pieces $V^i_k$ as in Section~\ref{unstable norm} where for each
$j$, the pieces $U^1_j$ and $U^2_j$ are paired up so that
$d_\Sigma (T^{-n}U^1_j, T^{-n}U^2_j) \leq C\lambda^{-n}\ve$. 
The introduction of the hole only increases the number of unpaired pieces $V^i_k$: if part of $T^{-n}W^1$
has fallen in the hole while the corresponding part of $T^{-n}W^2$ has not, then a piece $V^2_k \subset W^2$
is created.  We estimate the size of $V^2_k$ using assumption (H1).  

Suppose the part of $W^1$ 
corresponding to $V_k^2$ falls in the hole at time $\ell \leq n$.  If the boundary of the hole
at that point is strictly convex with curvature greater than $B$,  then
$|T^{-\ell}V_k^2| \leq C \sqrt{\ve}$ and so $|V_k^2| \leq C \sqrt{\ve}$ as well.
On the other hand, if the boundary of the hole is transverse to the stable cone, then the estimate
improves to $|T^{-\ell}V_k^2| \leq C \ve$. 
Notice also that there can be at most $L_n + P_n + 2$ pieces $V^i_k$.

Using this bound on the $V^i_k$, \eqref{eq:first unstable} becomes,
\begin{equation}
\label{eq:unstable hole}
\begin{split}
\sum_{i,k} \int_{T^{-n}V^i_k}h|DT^n|^{-1}J_WT^n\vf_i\circ T^n \, dm
& \leq C \| h \|_s \sum_{i,k} |V^i_k|^\alpha |\,|DT^n|^{-1} (J_WT^n)^{1-\alpha}|_{C^0} \\
& \leq C \ve^{\alpha/2} \|h\|_s (L_n+P_n) \lambda^{-n} \mu^{-\alpha n} .
\end{split}
\end{equation}
The estimates on the paired pieces $U^i_j$ do not change so putting
together equation~\eqref{eq:unstable hole} 
with \eqref{eq:second unstable} and \eqref{eq:third unstable}, and
using the fact that $\alpha/2 \geq \beta$, we have
\[
\| \Lp_H^n h \|_u \leq C \lambda^{-\beta n} D_n \|h\|_u  
+ C\|h\|_s (D_n + (L_n + P_n)\lambda^{-n}\mu^{-\alpha n}).
\]
This completes the proof of \eqref{eq:su hole}.
\end{proof}

\appendix
\section{Distortion Bounds}
\label{distortion}

The following are distortion bounds used in deriving the Lasota-Yorke estimates 
which are standard for uniformly hyperbolic $\C^2$ maps.
For any $n\in\N$ and $x, y\in K\in\mathcal K_n$ the following 
estimates hold.
\begin{equation}
\label{eq:distortion}
\begin{split}
\left| \frac{|DT^n(x)|}{|DT^n(y)|} -1 \right| & \; \; \leq \; \; C
\max\{d(x,y), d(T^nx,T^n,y)\}     \\
\left| \frac{J_WT^n(x)}{J_WT^n(y)} -1 \right| & \; \; \leq \; \; C
\max\{d(x,y), d(T^nx,T^n,y)\} 
\end{split}
\end{equation}
In particular, these bounds imply that $||DT^n|^{-1}|_{C^q(W_i)} \leq C ||DT^n|^{-1}|_{\C^0(W_i)} $ and 
similarly $|J_WT^n|_{C^q(W_i)} \leq C |J_WT^n|_{\C^0(W_i)}$ for any $0 \leq q \leq 1$.

Note that for $x \in T^{-n}W$, $|DT^n(x)| = C_{\theta}(x)J_WT^n(x)
J_uT^n(x)$ where $J_uT^n$ is the Jacobian  
of $T^n$ in the unstable direction and $C_{\theta}(x)$ is a number
which depends on the angle between the unstable  
direction and $T^{-n}W$ at the point $x$.  Since the family of
admissible leaves $W$ is uniformly transversal to the unstable
direction, there exists a constant $c_0>0$, independent of $W$, such
that $|C_{\theta}| \geq c_0$. Thus for all $n \geq 0$,
\begin{equation}
\label{eq:angle dist}
\left| |DT^n|^{-1} J_WT^n \right|_{\infty} \leq C \lambda^{-n}
\end{equation}
wherever $|DT^n|$ is defined.


\begin{thebibliography}{999999}
\footnotesize
\bibitem[Ba]{Bakh}  V.I. Bakhtin, {\em A direct method for constructing an invariant measure 
    on a hyperbolic attractor}.  (Russian) Izv. Ross. Akad. Nauk Ser. Mat. {\bf 56} (1992), 934-957;
    English transl., Russian Acad. Sci. Izv. Math. {\bf 41}:2 (1993), 207-227.
\bibitem[B1]{Bal} V. Baladi, {\em  Positive transfer operators and decay
    of correlations}, Advanced Series in Nonlinear Dynamics, {\bf 16},
    World Scientific (2000).
\bibitem[B2]{Ba1} V. Baladi, {\em Anisotropic Sobolev spaces and
    dynamical transfer operators: $\C^\infty$ foliations},  Algebraic and Topological Dynamics, Sergiy 
    Kolyada, Yuri Manin and Tom Ward, eds.  Contemporary Mathematics, Amer. Math. Society, (2005) 
    123-136.
\bibitem[BT]{BaT} V. Baladi, M. Tsujii, {\em Anisotropic H\"older 
    and Sobolev spaces for hyperbolic diffeomorphisms}, to appear in Ann. Inst. Fourier.
\bibitem[BY]{BalYo} V.Baladi, L.-S.Young, {\em On the spectra of randomly
    perturbed expanding maps}, Comm. Math. Phys., {\bf 156}:2 (1993),
    355-385; {\bf 166}:1 (1994), 219-220.
\bibitem[BC]{chernov bedem}  H. van den Bedem and N. Chernov, \emph{Expanding maps of an interval with holes}, 
    Ergod. Th. and Dynam. Sys. {\bf 22} (2002), 637-654.
\bibitem[BKL]{BKL} M. Blank, G. Keller, C. Liverani, {\em
    Ruelle-Perron-Frobenius spectrum for Anosov maps}, Nonlinearity, {\bf 15}:6 (2001), 1905-1973.
\bibitem[Bu]{buzzi} J. Buzzi,\emph{Absolutely continuous invariant probability measures for 
    arbitrary expanding piecewise $\mathbb{R}$-analytic mappings of the plane},
    Ergod. Th. and Dynam. Sys. {\bf 20}:3 (2000) 697-708.
\bibitem[BK]{buzzi keller} J. Buzzi and G. Keller, \emph{Zeta functions and transfer operators 
    for multidimensional piecewise affine and expanding maps}, Ergod. Th. and Dynam. Sys.  
    {\bf 21}:3  (2001), 689-716.
\bibitem[C]{cencova} N. N. \v{C}encova, \emph{A natural invariant
measure on Smale's horseshoe}, Soviet Math. Dokl. {\bf 23} (1981),
87-91.
\bibitem[CG]{chazottes}J.-R. Chazottes and S. Gouezel, \emph{On almost-sure versions of classical theorems for
dynamical systems}, to appear in Probability Theory and Related Fields.
\bibitem[Ch]{chernov 0} N. Chernov, \emph{Advanced statistical properties of
dispersing billiards}, J. Stat. Phys. {\bf 122} (2006), 1061-1094.
\bibitem[CD]{chernov dolgopyat} N. Chernov, D.Dolgopyat, \emph{Brownian Brownian Motion -
I}, to appear in Memoirs of AMS.
\bibitem[CM1]{chernov 1}  N. Chernov and R. Markarian, \emph{Ergodic
properties of Anosov maps with rectangular holes},
Bol. Soc. Bras. Mat. {\bf 28} (1997), 271-314. 
\bibitem[CM2]{chernov 2}  N. Chernov and R. Markarian, \emph{Anosov
maps with rectangular holes.  Nonergodic cases}, Bol. Soc. Bras. Mat. {\bf 28} (1997), 315-342.
\bibitem[CMT1]{chernov 3}  N. Chernov, R. Markarian and
S. Troubetskoy, \emph{Conditionally invariant measures for Anosov maps
with small holes}, Ergod. Th. and Dynam. Sys. {\bf 18} (1998),
1049-1073. 
\bibitem[CMT2]{chernov 4}  N. Chernov, R. Markarian and S. Troubetskoy,
    \emph{Invariant measures for Anosov maps with small holes},
Ergod. Th. and Dynam. Sys. {\bf 20} (2000), 1007-1044. 
\bibitem[CY]{chernov young}  N. Chernov and L.-S. Young, \emph{Decay of correlations for
    Lorentz gases and hard balls}, in Hard Ball Systems and the Lorentz Gas, D.Szasz, ed., 
    Enclyclopaedia of Mathematical Sciences {\bf 101},  Springer-Verlag:Berlin, 2000, 89-120. 
\bibitem[D1]{demers exp}  M. Demers, \emph{Markov Extensions for Dynamical Systems with Holes:  An Application to
    Expanding Maps of the Interval}, Israel J. of Math. {\bf 146} (2005), 189-221.
\bibitem[D2]{demers logistic}  M. Demers, \emph{Markov Extensions and Conditionally Invariant Measures for 
    Certain Logistic Maps with Small Holes}, Ergod. Th. and Dynam. Sys. {\bf 25}:4 (2005), 1139-1171.
\bibitem[DY]{demers young} M. Demers and L.-S. Young, \emph{Escape rates and natural conditionally
    invariant measures}, Nonlinearity, {\bf 19} (2006), 377-397.
\bibitem[GL]{liverani gouezel} S. Gou\"{e}zel and C. Liverani,
    \emph{Banach spaces adapted to Anosov systems},  Ergod.
    Th. and Dynam. Sys., {\bf 26}, 1, 189--217 (2006). 
\bibitem[HH]{hennion} H. Hennion, L.Hevr\'e, {\em Limit Theorems for
Markov chains and stochastic properties of dynamical systems by
quasi-compactness}, {\bf 1766}, Lectures Notes in Mathematics,
Springer-Verlag, Berlin, 2001.
\bibitem[K]{Kel1} G. Keller,
    {\em On the rate of convergence to equilibrium in one-dimensional systems}, 
    Comm. Math. Phys. {\bf 96} (1984), no. 2, 181--193.
\bibitem[KL]{KL} G. Keller, C. Liverani, {\em Stability of the spectrum for
    transfer operators}, Annali della Scuola Normale Superiore di Pisa,
    Scienze Fisiche e Matematiche, (4) {\bf XXVIII} (1999), 141-152.
\bibitem[LY]{LY} A. Lasota and J.A. Yorke, {\em On the existence of
    invariant measures for piecewise monotonic transformations},
    Trans. Amer. Math. Soc. {\bf 186} (1963), 481-488.
\bibitem[L1]{Li1}  C. Liverani, {\em  Decay of Correlations}, Annals of
    Mathematics {\bf 142} (1995), 239-301.
\bibitem[L2]{Li2} C. Liverani, {\em  Invariant measures and their
    properties. A functional analytic point of view}, Dynamical
    Systems. Part II: Topological Geometrical and Ergodic Properties of
    Dynamics. Pubblicazioni della Classe di Scienze, Scuola Normale
    Superiore, Pisa. Centro di Ricerca Matematica "Ennio De Giorgi" :
    Proceedings. Published by the Scuola Normale Superiore in Pisa
    (2004).
\bibitem[L3]{Li3} C. Liverani, {\em Fredholm determinants, Anosov maps and
    Ruelle resonances} , Discrete and Continuous Dynamical Systems,  {\bf 13}:5 (2005), 1203-1215.
\bibitem[LiM]{liverani maume} C. Liverani and V. Maume-Deschamps, \emph{Lasota-Yorke maps with holes: conditionally 
    invariant probability measures and invariant probability measures on the survivor set}, Annales de l'Institut Henri
    Poincar\'{e} Probability and Statistics, {\bf 39} (2003), 385-412. 
\bibitem[LW]{liverani wojtkowski} C. Liverani and M. Wojtkowski,
    \emph{Ergodicity in Hamiltonian Systems}, Dynamics Reported,  {\bf 4} (1995), 130-202.
\bibitem[LM]{lopes markarian} A. Lopes and R. Markarian, \emph{Open billiards: cantor sets, invariant
    and conditionally invariant probabilities}, SIAM J. Appl. Math. {\bf 56} (1996), 651-680.
\bibitem[PP1]{parry poll 1} W. Parry and M. Pollicott, \emph{An analogue of the prime number 
    theorem for closed orbits of Axiom A flows}, Annals of Math. {\bf 118}:3  (1983), 573-591.
\bibitem[PP2]{parry poll 2} W. Parry and M. Pollicott, \emph{Zeta functions and the periodic 
    orbit structure of hyperbolic dynamics}, Ast\'erisque No. 187-188, (1990), 268 pp.
\bibitem[P]{pesin} Ya.B.Pesin, \emph{Dynamical systems with generalized
    hyperbolic attractors: hyperbolic, ergodic and topological
    properties}, Ergod. Th. and Dynam. Sys. {\bf 12} (1992), 123-151.
\bibitem[Ru1]{ruelle 1} D. Ruelle, \emph{Locating resonances for Axiom A 
    dynamical systems}, J. Stat. Phys.  {\bf 44}:3-4  (1986), 281-292.
\bibitem[Ru2]{ruelle 2} D. Ruelle, \emph{Resonances for Axiom $A$ flows}, 
    J. Differential Geom. {\bf 25}:1  (1987), 99-116.
\bibitem[R1]{Rugh1} H.H. Rugh, {\em The correlation spectrum for hyperbolic
    analytic maps},  Nonlinearity {\bf 5}:6  (1992),  1237-1263.
\bibitem[R2]{Rugh2} H.H. Rugh, {\em Fredholm determinants for
    real-analytic hyperbolic diffeomorphisms of surfaces}.  XIth
    International Congress of Mathematical Physics (Paris, 1994),
    297--303, Internat. Press, Cambridge, MA, 1995.
\bibitem[R3]{Rugh3} H.H. Rugh, {\em Generalized Fredholm determinants
    and Selberg zeta functions for Axiom A dynamical systems}.
    Ergod. Th. and Dynam. Sys.  {\bf 16}:4  (1996), 805-819.
\bibitem[S]{saussol}  B. Saussol, \emph{Absolutely continuous invariant measures for multidimensional 
    expanding maps}, Israel J. Math. {\bf 116} (2000), 223-248.
\bibitem[T1]{tsujii 1} M. Tsujii, \emph{Absolutely continuous invariant measures for piecewise 
    real-analytic expanding maps on the plane}, Comm. Math. Phys. {\bf 208}:3 (2000), 605-622.
\bibitem[T2]{tsujii 2} M. Tsujii, \emph{Absolutely continuous invariant measures for 
    expanding piecewise linear maps}, Invent. Math. {\bf 143}:2 (2001), 349-373.
\bibitem[Y]{young} L.-S. Young, \emph{Statistical properties of
    dynamical systems with some hyperbolicity}, Annals of Math. {\bf 147}:3 (1998), 585-650.
\end{thebibliography}
\end{document}